\documentclass[10pt,a4paper,oneside]{article}

\usepackage[titletoc]{appendix}

\usepackage{amsmath,amssymb}

\usepackage{extarrows}
\usepackage{pdflscape}

\usepackage{tocloft}

\DeclareMathOperator{\Cone}{Cone}

\DeclareMathOperator{\fchar}{char}
\DeclareMathOperator{\Fil}{Fil}
\DeclareMathOperator{\Gal}{Gal}
\DeclareMathOperator{\Hom}{Hom}

\DeclareMathOperator{\Isom}{Isom}
\DeclareMathOperator{\ord}{ord}

\DeclareMathOperator{\rk}{rk}
\DeclareMathOperator{\Spec}{Spec}

\newcommand{\CC}{\mathbb{C}}
\newcommand{\FF}{\mathbb{F}}

\newcommand{\QQ}{\mathbb{Q}}
\newcommand{\RR}{\mathbb{R}}
\newcommand{\ZZ}{\mathbb{Z}}

\renewcommand{\AA}{\mathbb{A}}
\newcommand{\PP}{\mathbb{P}}

\DeclareMathOperator{\Gr}{Gr}

\newcommand{\Parf}{\mathcal{P}\!\text{\it arf}}

\renewcommand{\det}{\operatorname{det}}

\renewcommand{\Re}{\operatorname{Re}}
\renewcommand{\emptyset}{\varnothing}

\newcommand{\ar}{\text{\it ar}}
\newcommand{\BM}{\text{\it BM}}

\newcommand{\DB}{{\mathcal{D}\text{-}\mathcal{B}}}

\newcommand{\dR}{\text{\it dR}}
\newcommand{\et}{\text{\it \'{e}t}}
\newcommand{\fg}{\text{\it fg}}
\newcommand{\is}{\text{\it is}}
\newcommand{\red}{\text{\it red}}

\newcommand{\Wc}{\text{\it W,c}}

\newcommand{\dfn}{\mathrel{\mathop:}=}

\newcommand{\iHom}{\underline{\Hom}}
\newcommand{\RHom}{R\!\Hom}

\usepackage{tikz-cd}
\usetikzlibrary{arrows}
\usetikzlibrary{calc}
\usetikzlibrary{babel}
\usetikzlibrary{decorations.pathreplacing}
\usetikzlibrary{patterns}
\newcommand{\tikzpb}{\ar[phantom,pos=0.2]{dr}{\text{\large$\lrcorner$}}}

\usepackage{hyperref}

\hypersetup{
    colorlinks,
    linkcolor={red!60!black},
    citecolor={blue!60!black},
    urlcolor={blue!80!black}
}

\usepackage{amsthm}

\newtheoremstyle{myplain}
{\topsep}   % ABOVESPACE
{\topsep}   % BELOWSPACE
{\itshape}  % BODYFONT
{0pt}       % INDENT (empty value is the same as 0pt)
{\scshape} % HEADFONT
{.}         % HEADPUNCT
{5pt plus 1pt minus 1pt} % HEADSPACE
{}   % CUSTOM-HEAD-SPEC

\theoremstyle{myplain}

\newtheorem*{maintheorem*}{Main theorem}
\newtheorem*{thetheorem*}{Theorem}
\newtheorem*{proposition*}{Proposition}
\newtheorem{theorem}{Theorem}[section]
\newtheorem{proposition}[theorem]{Proposition}
\newtheorem{lemma}[theorem]{Lemma}
\newtheorem{corollary}[theorem]{Corollary}

\newtheoremstyle{mydefinition}
{\topsep}   % ABOVESPACE
{\topsep}   % BELOWSPACE
{}  % BODYFONT
{0pt}       % INDENT (empty value is the same as 0pt)
{\scshape} % HEADFONT
{.}         % HEADPUNCT
{5pt plus 1pt minus 1pt} % HEADSPACE
{}   % CUSTOM-HEAD-SPEC

\theoremstyle{mydefinition}
\newtheorem*{conjecture*}{Conjecture}
\newtheorem{definition}[theorem]{Definition}
\newtheorem{conjecture}[theorem]{Conjecture}
\newtheorem{remark}[theorem]{Remark}
\newtheorem{example}[theorem]{Example}

\title{Weil-\'{e}tale cohomology and zeta-values of arithmetic schemes at negative integers}
\author{Alexey Beshenov}
\date{}

\AtEndDocument{%
  \par
  \medskip
  \begin{tabular}{@{}l@{}}%
    \\
    Alexey Beshenov \\
    E-mail: \texttt{cadadr@gmail.com} \\
    URL: \url{https://cadadr.org/}
  \end{tabular}}

\numberwithin{equation}{section}

\begin{document}

\maketitle

\begin{abstract}
  Following the ideas of Flach and Morin \cite{Flach-Morin-2018}, we state a
  conjecture in terms of Weil-\'{e}tale cohomology for the vanishing order and
  special value of the zeta function $\zeta (X,s)$ at $s = n < 0$, where $X$ is
  a separated scheme of finite type over $\Spec \ZZ$. We prove that the
  conjecture is compatible with closed-open decompositions of schemes and with
  affine bundles, and consequently, that it holds for cellular schemes over
  certain one-dimensional bases.

  This is a continuation of \cite{Beshenov-Weil-etale-1}, which gives a
  construction of Weil-\'{e}tale cohomology for $n<0$ under the mentioned
  assumptions on $X$.
\end{abstract}

% \tableofcontents

%%%%%%%%%%%%%%%%%%%%%%%%%%%%%%%%%%%%%%%%%%%%%%%%%%%%%%%%%%%%%%%%%%%%%%%%%%%%%%%%

\section{Introduction}

Let $X$ be an \textbf{arithmetic scheme}, by which we mean in this paper that it
is a separated scheme of finite type $X \to \Spec \ZZ$. Then the corresponding
\textbf{zeta function} is defined by
\begin{equation}
  \label{eqn:Euler-product-for-zeta}
  \zeta (X,s) = \prod_{\substack{x \in X \\ \text{closed pt.}}}
  \frac{1}{1 - N (x)^{-s}}.
\end{equation}
Here, for a closed point $x \in X$, the norm
$$N (x) = |\kappa (x)| = |\mathcal{O}_{X,x}/\mathfrak{m}_{X,x}|$$
is the size of the corresponding residue field. The product converges for
$\Re s > \dim X$, and conjecturally admits a meromorphic continuation to the
whole complex plane. Basic facts and conjectures about zeta functions of schemes
can be found in \cite{Serre-1965}.

Of particular interest are the so-called special values of $\zeta (X,s)$ at
integers $s = n \in \ZZ$, also known as the \textbf{zeta-values} of $X$.
To define these, we assume that $\zeta (X,s)$ admits a meromorphic continuation
around $s = n$. We denote by
$$d_n = \ord_{s=n} \zeta (X,s)$$
the \textbf{vanishing order} of $\zeta (X,s)$ at $s = n$. That is, $d_n > 0$
(resp. $d_n < 0$) if $\zeta (X,s)$ has a zero (resp. pole) of order $|d_n|$ at
$s = n$.

The \textbf{special value} of $\zeta (X,s)$ at $s = n$ is defined as the leading
nonzero coefficient of the Taylor expansion:
$$\zeta^* (X,n) = \lim_{s \to n} (s - n)^{-d_n}\,\zeta (X,s).$$

Early on, Lichtenbaum conjectured that both numbers $\ord_{s = n} \zeta (X,s)$
and $\zeta^* (X,n)$ should have a cohomological interpretation related to the
\'{e}tale motivic cohomology of $X$ (see e.g. \cite{Lichtenbaum-1984} for
varieties over finite fields).

This is made precise in Lichtenbaum's Weil-\'{e}tale program. It suggests the
existence of \textbf{Weil-\'{e}tale cohomology}, which is a suitable
modification of motivic cohomology that encodes the information about the
vanishing order and the special value of $\zeta (X,s)$ at $s = n$.
For Lichtenbaum's recent work on this topic, we refer the reader to
\cite{Lichtenbaum-2005,Lichtenbaum-2009-Euler-char,Lichtenbaum-2009-number-rings,Lichtenbaum-2021}.

The case of varieties over finite fields $X/\FF_q$ is now well understood thanks
to the work of Geisser
\cite{Geisser-2004,Geisser-2006,Geisser-2010-arithmetic-homology}.

Flach and Morin considered the case of proper, regular arithmetic schemes
$X$. In \cite{Flach-Morin-2012} they have studied the corresponding
Weil-\'{e}tale topos. Later, in \cite{Morin-2014} Morin gave an explicit
construction of Weil-\'{e}tale cohomology groups $H^i_\Wc (X, \ZZ)$ for a proper
and regular arithmetic scheme $X$. This construction was further generalized by
Flach and Morin in \cite{Flach-Morin-2018} to groups $H^i_\Wc (X, \ZZ(n))$ with
weights $n \in \ZZ$, again for a proper and regular $X$.

Motivated by the work of Flach and Morin, the author constructed in
\cite{Beshenov-Weil-etale-1} Weil-\'{e}tale cohomology groups
$H^i_\Wc (X, \ZZ (n))$ for any arithmetic scheme $X$ (removing the assumption
that $X$ is proper or regular) and strictly negative weights $n < 0$.
The construction is based on the following assumption.

\begin{conjecture*}
  $\mathbf{L}^c (X_\et,n)$: given an arithmetic scheme $X$ and $n < 0$, the
  cohomology groups $H^i (X_\et, \ZZ^c (n))$ are finitely generated for all
  $i \in \ZZ$.
\end{conjecture*}

For the known cases, see \cite[\S 8]{Beshenov-Weil-etale-1}. Under this
conjecture, we constructed in \cite[\S 7]{Beshenov-Weil-etale-1} perfect
complexes of abelian groups $R\Gamma_\Wc (X, \ZZ(n))$ and the corresponding
cohomology groups
$$H^i_\Wc (X, \ZZ(n)) \dfn H^i (R\Gamma_\Wc (X, \ZZ(n))).$$

This text is a continuation of \cite{Beshenov-Weil-etale-1} and investigates the
conjectural relation of our Weil-\'{e}tale cohomology to the special value of
$\zeta (X,s)$ at $s = n < 0$.  Specifically, we make the following conjectures.

\begin{enumerate}
\item[1)] \textbf{Conjecture}~$\mathbf{VO} (X,n)$
  (see \S\ref{sec:vanishing-order-conjecture}):
  \emph{the vanishing order is given by the weighted alternating sum of ranks}
  \[ \ord_{s=n} \zeta (X,s) =
    \sum_{i\in \ZZ} (-1)^i \cdot i \cdot \rk_\ZZ H_\Wc^i (X, \ZZ(n)). \]

\item[2)] A consequence of \textbf{Conjecture}~$\mathbf{B} (X,n)$
  (see \S\ref{sec:regulator} and Lemma~\ref{lemma:smile-theta}):
  \emph{after tensoring the cohomology groups $H_\Wc^i (X, \ZZ (n))$ with $\RR$,
    we obtain a long exact sequence of finite dimensional real vector spaces}
  \[ \cdots \to H_\Wc^{i-1} (X, \RR (n)) \xrightarrow{\smile\theta}
    H_\Wc^i (X, \RR (n)) \xrightarrow{\smile\theta}
    H_\Wc^{i+1} (X, \RR (n)) \to \cdots \]

  It follows that there is a canonical isomorphism
  \[ \lambda\colon \RR \xrightarrow{\cong}
    (\det_\ZZ R\Gamma_\Wc (X, \ZZ (n))) \otimes \RR. \]
  Here $\det_\ZZ R\Gamma_\Wc (X, \ZZ (n))$ is the determinant of the
  perfect complex of abelian groups $R\Gamma_\Wc (X, \ZZ (n))$, in the sense of
  Knudsen and Mumford \cite{Knudsen-Mumford-1976}. In particular,
  $\det_\ZZ R\Gamma_\Wc (X, \ZZ (n))$ is a free $\ZZ$-module of rank
  $1$. For the convenience of the reader, we give a brief overview of
  determinants in Appendix~\ref{app:determinants}.

\item[3)] \textbf{Conjecture}~$\mathbf{C} (X,n)$
  (see \S\ref{sec:special-value-conjecture}):
  \emph{the special value is determined up to sign by}
  \[ \lambda (\zeta^* (X, n)^{-1}) \cdot \ZZ =
    \det_\ZZ R\Gamma_\Wc (X, \ZZ (n)). \]
\end{enumerate}

If $X$ is proper and regular, then our construction of
$R\Gamma_\Wc (X, \ZZ (n))$ and the above conjectures agree with those of Flach
and Morin from \cite{Flach-Morin-2018}. Apart from removing the assumption that
$X$ is proper and regular, a novelty of this work is that we prove the
compatibility of the conjectures with operations on schemes, in particular with
closed-open decompositions $Z \not\hookrightarrow X \hookleftarrow U$, where
$Z \subset X$ is a closed subscheme and $U = X\setminus Z$ is the open
complement, and with affine bundles $\AA_X^r = \AA_\ZZ^r \times X$ (see
Proposition~\ref{prop:compatibility-of-VO(X,n)} and
Theorem~\ref{thm:compatibility-of-C(X,n)}). This gives a machinery for starting
from the cases of schemes for which the conjectures are known and constructing
new schemes for which the conjectures also hold.  As an application, we prove in
\S\ref{sec:unconditional-results} the following result.

\begin{maintheorem*}
  Let $B$ be a one-dimensional arithmetic scheme, such that each of the generic
  points $\eta \in B$ satisfies one of the following properties:
  \begin{enumerate}
  \item[a)] $\fchar \kappa (\eta) = p > 0$;

  \item[b)] $\fchar \kappa (\eta) = 0$, and $\kappa (\eta)/\QQ$ is an abelian
    number field.
  \end{enumerate}
  If $X$ is a $B$-cellular arithmetic scheme with smooth quasi-projective fiber
  $X_{\red,\CC}$, then Conjectures~$\mathbf{VO} (X,n)$ and
  $\mathbf{C} (X,n)$ hold unconditionally for any $n < 0$.
\end{maintheorem*}

In fact, this result is established for a larger class of arithmetic schemes
$\mathcal{C} (\ZZ)$; we refer to \S\ref{sec:unconditional-results} for more
details.

\subsection*{Outline of the paper}

In \S\ref{sec:regulator} we define the regulator morphism, based on the
construction of Kerr, Lewis, and M\"{u}ller-Stach
\cite{Kerr-Lewis-Muller-Stach-2006}, and state the associated
Conjecture~$\mathbf{B} (X,n)$.

Then \S\ref{sec:vanishing-order-conjecture} is devoted to
Conjecture~$\mathbf{VO} (X,n)$ about the vanishing order. We also explain why it
is consistent with a conjecture of Soul\'{e}, and with the vanishing order
arising from the expected functional equation.

In \S\ref{sec:special-value-conjecture} we state Conjecture~$\mathbf{C} (X,n)$
about the special value.

We explain in \S\ref{sec:finite-fields} that if $X$ is a variety over a finite
field, then Conjecture~$\mathbf{C} (X,n)$ is consistent with the conjectures
considered by Geisser in
\cite{Geisser-2004,Geisser-2006,Geisser-2010-arithmetic-homology}, and it
follows from Conjecture~$\mathbf{L}^c (X_\et,n)$.

Then we prove in \S\ref{sec:compatibility-with-operations} that Conjectures
$\mathbf{VO} (X,n)$ and $\mathbf{C} (X,n)$ are compatible with basic operations
on schemes: disjoint unions, closed-open decompositions, and affine
bundles. Using these results, we conclude in \S\ref{sec:unconditional-results}
with a class of schemes for which the conjectures hold unconditionally.

For the convenience of the reader, Appendix~\ref{app:determinants} gives a brief
overview of basic definitions and facts related to the determinants of
complexes.

\subsection*{Notation}

In this paper, $X$ always denotes an \textbf{arithmetic scheme} (separated, of
finite type over $\Spec \ZZ$), and $n$ is always a strictly negative integer.

We denote by
\[
  R\Gamma_\fg (X, \ZZ (n))
  \quad\text{and}\quad
  R\Gamma_\Wc (X, \ZZ (n))
\]
the complexes of abelian groups constructed in \cite{Beshenov-Weil-etale-1}
under Conjecture~$\mathbf{L}^c (X_\et,n)$. We set
\begin{align*}
  H^i_\fg (X, \ZZ (n)) & \dfn H^i (R\Gamma_\fg (X, \ZZ (n))), \\
  H^i_\Wc (X, \ZZ (n)) & \dfn H^i (R\Gamma_\Wc (X, \ZZ (n))).
\end{align*}

By \cite[Proposition 5.5 and 7.12]{Beshenov-Weil-etale-1}, these
cohomology groups are finitely generated, assuming
Conjecture~$\mathbf{L}^c (X_\et,n)$; moreover, the groups $H^i_\Wc (X, \ZZ(n))$
are bounded, and $H^i_\fg (X, \ZZ (n))$ are bounded from below and finite
$2$-torsion for $i \gg 0$.

Briefly, the construction fits in the following diagram of distinguished
triangles in the derived category $\mathbf{D} (\ZZ)$:
\[ \begin{tikzcd}[column sep=1.5em]
    &[-3em] \RHom (R\Gamma (X_\et, \ZZ^c (n)), \QQ [-2]) \ar{d}{\alpha_{X,n}} \ar{r} &[-2.5em] 0 \ar{d} \\
    & R\Gamma_c (X_\et, \ZZ(n)) \ar{d}\ar{r}{u_\infty^*} & R\Gamma_c (G_\RR, X (\CC), \ZZ(n))\ar{d}{id} \\
    R\Gamma_\Wc (X, \ZZ (n)) \ar{r} & R\Gamma_\fg (X, \ZZ(n)) \ar[dashed]{r}{i_\infty^*}\ar{d} & R\Gamma_c (G_\RR, X (\CC), \ZZ(n)) \ar{r} \ar{d} & {[1]} \\
    & \RHom (R\Gamma (X_\et, \ZZ^c (n)), \QQ [-1]) \ar{r} & 0
\end{tikzcd} \]
For more details, see \cite{Beshenov-Weil-etale-1}.
For real coefficients, we set
\begin{align*}
  R\Gamma_\fg (X, \RR (n)) & \dfn R\Gamma_\fg (X, \ZZ (n)) \otimes \RR, \\
  R\Gamma_\Wc (X, \RR (n)) & \dfn R\Gamma_\Wc (X, \ZZ (n)) \otimes \RR.
\end{align*}
Accordingly,
\begin{align*}
  H^i_\fg (X, \RR (n)) & \dfn H^i (R\Gamma_\fg (X, \RR (n))) = H^i_\fg (X, \ZZ (n)) \otimes \RR, \\
  H^i_\Wc (X, \RR (n)) & \dfn H^i (R\Gamma_\Wc (X, \RR (n))) = H^i_\Wc (X, \ZZ (n)) \otimes \RR.
\end{align*}

By $X (\CC)$ we denote the space of complex points of $X$ with the usual
analytic topology. It carries a natural action of $G_\RR = \Gal (\CC/\RR)$ via
the complex conjugation. For a subring $A \subseteq \RR$ we denote by $A (n)$
the $G_\RR$-module $(2\pi i)^n\,A$, and also the corresponding constant
$G_\RR$-equivariant sheaf on $X (\CC)$.

We denote by $R\Gamma_c (X (\CC), A (n))$ the cohomology with compact support
with $A (n)$-coefficients, and its $G_\RR$-equivariant version is defined by
$$R\Gamma_c (G_\RR, X (\CC), A (n)) \dfn R\Gamma (G_\RR, R\Gamma_c (X (\CC), A (n)))$$
For real coefficients, we have
$$H_c^i (G_\RR, X (\CC), \RR (n)) = H^i_c (X (\CC), \RR (n))^{G_\RR},$$
where the $G_\RR$-action on $H^i_c (X (\CC), \RR (n))$ naturally comes from the
corresponding action on $X (\CC)$ and $\RR (n)$.

\textbf{Borel--Moore homology} is defined as the dual to cohomology with compact
support. We are interested in the real coefficients:
\begin{align*}
  R\Gamma_\BM (X (\CC), \RR (n)) & \dfn
  \RHom (R\Gamma_c (X (\CC), \RR (n)), \RR), \\
  R\Gamma_\BM (G_\RR, X (\CC), \RR (n)) & \dfn
  \RHom (R\Gamma_c (G_\RR, X (\CC), \RR (n)), \RR).
\end{align*}

\subsection*{Acknowledgments}

Parts of this work are based on my doctoral thesis, which I wrote under the
supervision of Baptiste Morin (Universit\'{e} de Bordeaux) and Bas Edixhoven
(Universiteit Leiden). I am very grateful to them for their support in working
on this project. I am also indebted to Matthias Flach, as the ideas for this
work came from \cite{Flach-Morin-2018}. I thank Stephen Lichtenbaum and Niranjan
Ramachandran who kindly agreed to act as reviewers for my thesis and provided me
with many useful comments and suggestions. Finally, I thank Jos\'{e} Jaime
Hern\'{a}ndez Castillo, Diosel L\'{o}pez Cruz, and Maxim Mornev for several
fruitful discussions.

This paper was edited while I visited the Center for Research in Mathematics
(CIMAT), Guanajuato. I personally thank Pedro Luis del \'{A}ngel and
Xavier G\'{o}mez Mont for their hospitality.

%%%%%%%%%%%%%%%%%%%%%%%%%%%%%%%%%%%%%%%%%%%%%%%%%%%%%%%%%%%%%%%%%%%%%%%%%%%%%%%%

\section{Regulator morphism and Conjecture~$\mathbf{B} (X,n)$}
\label{sec:regulator}

In order to formulate the special value conjecture, we need a regulator morphism
from motivic cohomology to Deligne(--Beilinson) (co)homology. Such regulators
were originally introduced by Bloch in \cite{Bloch-1986-Lefschetz}, and here we
use the construction of Kerr, Lewis, and M\"{u}ller-Stach
\cite{Kerr-Lewis-Muller-Stach-2006}, which works at the level of complexes.
We will simply call it ``the KLM regulator.'' It works under the assumption that
$X_{\red,\CC}$ is a smooth quasi-projective variety.

For simplicity, in this section we assume that $X$ is reduced (motivic
cohomology does not distinguish between $X$ and $X_\red$), and that $X_\CC$ is
connected of dimension $d_\CC$ (otherwise, the arguments below can be applied to
each connected component). We fix a compactification by a normal crossing
divisor
\[ \begin{tikzcd}
    X_\CC \ar[right hook->]{r}{j} & \overline{X}_\CC & D \ar[left hook->]{l}
  \end{tikzcd} \]
The KLM regulator has the form of a morphism in the derived category
\begin{equation}
  \label{eqn:KLM-morphism-1}
  z^p (X_\CC, -\bullet) \otimes \QQ \to
  {}' C_\mathcal{D}^{2p - 2d_\CC + \bullet} (\overline{X}_\CC, D, \QQ (p-d_\CC)).
\end{equation}

Here $z^p (X_\CC, -\bullet)$ denotes the Bloch's cycle complex
\cite{Bloch-1986}. To define it, consider the algebraic simplex
$\Delta_\CC^i = \Spec \CC [t_0,\ldots,t_i]/(1 - \sum_j t_j)$.
Then, $z^p (X_\CC, i)$ is freely generated by algebraic cycles
$Z \subset X_\CC \times_{\Spec \CC} \Delta_\CC^i$ of codimension $p$ which
intersect the faces properly. It is more convenient for us to work with
$$z_{d_\CC - p} (X_\CC, i) = z^p (X_\CC, i),$$
generated by cycles $Z \subset X_\CC \times_{\Spec \CC} \Delta_\CC^i$ of
dimension $p+i$.

The complex ${}' C_\mathcal{D}^{\bullet} (\overline{X}_\CC, D, \QQ (k))$ on the
right-hand side of \eqref{eqn:KLM-morphism-1} computes Deligne(--Beilinson)
homology, as defined by Jannsen \cite{Jannsen-1988}. If we take
$p = d_\CC + 1 - n$, tensor it with $\RR$ and shift it by $2n$, we obtain
\begin{equation}
  \label{eqn:KLM-morphism-2}
  z_{n-1} (X_\CC, -\bullet) \otimes \RR [2n] \to
  {}' C_\mathcal{D}^{2 + \bullet} (\overline{X}_\CC, D, \RR (1-n)).
\end{equation}

\begin{remark}
  Some comments are in order.

  \begin{enumerate}
  \item Originally, the KLM regulator is defined using a cubical version of
    cycle complexes, but these are quasi-isomorphic to the usual simplicial
    cycle complexes by \cite{Levine-1994}, so we make no distinction here.
    For an explicit simplicial version of the KLM regulator, see
    \cite{Kerr-Lewis-Lopatto-2018}.

  \item The KLM regulator is defined as a true morphism of complexes (not just a
    morphism in the derived category) on a subcomplex
    $z^r_\RR (X_\CC, \bullet) \subset z^r (X_\CC, \bullet)$. This inclusion
    becomes a quasi-isomorphism if we pass to rational coefficients. In the
    original paper \cite{Kerr-Lewis-Muller-Stach-2006} this is stated without
    tensoring with $\QQ$ , but the omission is acknowledged later in
    \cite{Kerr-Lewis-2007}. For our purposes, it suffices to have a regulator
    with coefficients in $\RR$.

  \item The case of a smooth quasi-projective $X_\CC$, where one must consider a
    compactification by a normal crossing divisor as above, is treated in
    \cite[\S 5.9]{Kerr-Lewis-Muller-Stach-2006}.
  \end{enumerate}
\end{remark}

Now we make a small digression to identify the right-hand side of
\eqref{eqn:KLM-morphism-2}. Under our assumption that $n < 0$, Deligne
homology is equivalent to Borel--Moore homology.

\begin{lemma}
  For any $n < 0$ there is a quasi-isomorphism
  \begin{multline*}
    {}' C^\bullet_\mathcal{D} (\overline{X_\CC}, D, \RR (1-n)) \cong
    R\Gamma_\BM (X (\CC), \RR (n)) [-1] \\
    \dfn \RHom (R\Gamma_c (X (\CC), \RR (n)), \RR) [-1].
  \end{multline*}
  Moreover, it respects the natural actions of $G_\RR$ on both complexes.

  \begin{proof}
    From the proof of \cite[Theorem~1.15]{Jannsen-1988}, for any $k \in \ZZ$ we
    have a quasi-isomorphism
    \begin{equation}
      \label{eqn:Jannsen-Theorem-1.15}
      {}' C^\bullet_\mathcal{D} (\overline{X_\CC}, D, \RR (k)) \cong
      R\Gamma (\overline{X} (\CC), \RR (k + d_\CC)_{\DB, (\overline{X}_\CC,X_\CC)}) [2d_\CC],
    \end{equation}
    where
    \[ \RR (k + d_\CC)_{\DB, (\overline{X}_\CC,X_\CC)} =
      \Cone \left(\begin{array}{c}
        R j_* \RR (k + d_\CC) \\
        \oplus \\
        \Omega^{\geqslant k + d_\CC}_{\overline{X} (\CC)} (\log D)
      \end{array}
      \xrightarrow{\epsilon - \iota}
      R j_* \Omega_{X (\CC)}^\bullet \right) [-1] \]
    is the sheaf whose hypercohomology on $\overline{X} (\CC)$ gives
    Deligne--Beilinson cohomology (see \cite{Esnault-Viehweg-1988} for
    more details).

    Here $\Omega^\bullet_{\overline{X} (\CC)}$ denotes the usual de Rham complex
    of holomorphic differential forms, and
    $\Omega^\bullet_{\overline{X} (\CC)} (\log D)$ is the complex of forms with
    at most logarithmic poles along $D (\CC)$.
    The latter complex is filtered by subcomplexes
    $\Omega^{\geqslant \bullet}_{\overline{X} (\CC)} (\log D)$.
    The morphism
    $\epsilon\colon R j_* \RR (k) \to R j_* \Omega^\bullet_{X (\CC)}$ is induced
    by the canonical morphism of sheaves $\RR (k) \to \mathcal{O}_{X (\CC)}$,
    and $\iota$ is induced by the natural inclusion
    $\Omega^\bullet_{\overline{X} (\CC)} (\log D) \xrightarrow{\cong} j_*
    \Omega_{X (\CC)}^\bullet = R j_* \Omega_{X (\CC)}^\bullet$, which is a
    quasi-isomorphism of filtered complexes.

    We are interested in the case of $k > 0$ when the part
    $\Omega^{\geqslant k + d_\CC}_{\overline{X} (\CC)} (\log D)$ vanishes, and
    we obtain
    \begin{align}
      \notag \RR (k + d_\CC)_{\DB, (\overline{X}_\CC,X_\CC)} & \cong
                                     R j_* \Cone \Bigl(\RR (k + d_\CC)
                                     \xrightarrow{\epsilon}
                                     \Omega_{X (\CC)}^\bullet \Bigr) [-1] \\
      \notag & \cong R j_* \Bigl(\RR (k + d_\CC) \xrightarrow{\epsilon}
               \Omega_{X (\CC)}^\bullet [-1] \Bigr) \\
      \label{eqn:deligne-homology-1} & \cong R j_* \Bigl(\RR (k + d_\CC) \to \CC [-1] \Bigr) \\
      \label{eqn:deligne-homology-2} & \cong R j_* \RR (k + d_\CC - 1) [-1]
    \end{align}
    Here \eqref{eqn:deligne-homology-1} comes from the Poincar\'{e} lemma
    $\CC \cong \Omega_{X (\CC)}^\bullet$ and \eqref{eqn:deligne-homology-2}
    from the short exact sequence of $G_\RR$-modules
    $\RR (k + d_\CC) \rightarrowtail \CC \twoheadrightarrow \RR (k + d_\CC - 1)$.

    Returning to \eqref{eqn:Jannsen-Theorem-1.15} for $k = 1-n$, we find that
    \begin{align*}
      {}' C^\bullet_\mathcal{D} (\overline{X_\CC}, D, \RR (1-n)) & \cong
      R\Gamma (X (\CC), \RR (d_\CC - n)) [2d_\CC-1] \\
      & \cong \RHom (R\Gamma_c (X (\CC), \RR (n)), \RR) [-1].
    \end{align*}
    Here the final isomorphism is Poincar\'{e} duality.
    All the above is $G_\RR$-equivariant.
  \end{proof}
\end{lemma}

Returning now to \eqref{eqn:KLM-morphism-2}, the previous lemma allows us to
reinterpret the KLM regulator as
\begin{equation}
  \label{eqn:KLM-morphism-3}
  z_{n-1} (X_\CC, -\bullet) \otimes \RR [2n] \to
  R\Gamma_\BM (X (\CC), \RR (n)), \RR) [1].
\end{equation}

We have
\begin{multline}
  \label{eqn:KLM-morphism-4}
  z_{n-1} (X_\CC, -\bullet) \otimes \RR [2n] =
  z_{n-1} (X_\CC, -\bullet) \otimes \RR [2n-2] [2] \\
  = \Gamma (X_{\CC,\et}, \RR^c (n-1)) [2],
\end{multline}
where the complex of sheaves $\RR^c (p)$ is defined by
$U \rightsquigarrow z_p (U, -\bullet) \otimes \RR [2p]$.
By \'{e}tale cohomological descent \cite[Theorem~3.1]{Geisser-2010},
\begin{equation}
  \label{eqn:KLM-morphism-5}
  \Gamma (X_{\CC,\et}, \RR^c (n-1)) \cong R\Gamma (X_{\CC,\et}, \RR^c (n-1)).
\end{equation}
(We note that \cite[Theorem~3.1]{Geisser-2010} holds unconditionally, since the
Beilinson--Lichtenbaum conjecture follows from the Bloch--Kato conjecture, which
is now a theorem; see also \cite{Geisser-2004-Dedekind} where the consequences
of Bloch--Kato for motivic cohomology are deduced.)

Finally, the base change from $X$ to $X_\CC$ naturally maps cycles
$Z \subset X \times \Delta_\ZZ^i$ of dimension $n$ to cycles in
$X_\CC \times_{\Spec \CC} \Delta_\CC^i$ of dimension $n-1$, so that there is a
morphism
\begin{equation}
  \label{eqn:KLM-morphism-6}
  R\Gamma (X_\et, \RR^c (n)) \to R\Gamma (X_{\CC,\et}, \RR^c (n-1)) [2].
\end{equation}

\begin{remark}
  Assuming that $X$ is flat and has pure Krull dimension $d$, we have
  $\RR^c (n)^X = \RR (d-n)^X [2d]$, where $\RR (\bullet)$ is the usual cycle
  complex defined by $z^n (\text{\textvisiblespace}, -\bullet) [-2n]$.
  Similarly, $\RR^c (n)^{X_\CC} = \RR (d_\CC-n)^{X_\CC} [2d_\CC]$, with
  $d_\CC = d - 1$. With this renumbering, the morphism
  \eqref{eqn:KLM-morphism-6} becomes
  $$R\Gamma (X_\et, \RR (d-n)) [2d] \to R\Gamma (X_{\CC,\et}, \RR (d-n)) [2d].$$
  This probably looks more natural, but we make no additional assumptions about
  $X$ and work exclusively with complexes $A^c (\bullet)$ defined in terms of
  dimension of algebraic cycles, rather than $A (\bullet)$ defined in terms of
  codimension.
\end{remark}

\begin{definition}
  Given an arithmetic scheme $X$ with smooth quasi-projective $X_\CC$ and
  $n < 0$, consider the composition of morphisms
  \begin{multline*}
    R\Gamma (X_\et, \RR^c (n)) \xrightarrow{\text{\eqref{eqn:KLM-morphism-6}}}
    R\Gamma (X_{\CC,\et}, \RR^c (n-1)) [2] \stackrel{\text{\eqref{eqn:KLM-morphism-5}}}{\cong}
    \Gamma (X_{\CC,\et}, \RR^c (n-1)) [2] \\
    \stackrel{\text{\eqref{eqn:KLM-morphism-4}}}{=}
    z_{n-1} (X_\CC, -\bullet)_\RR [2n] \xrightarrow{\text{\eqref{eqn:KLM-morphism-3}}}
    R\Gamma_\BM (X (\CC), \RR (n)), \RR) [1].
  \end{multline*}
  Moreover, we take the $G_\RR$-invariants, which gives us the
  \textbf{(\'{e}tale) regulator}
  \[ Reg_{X,n}\colon R\Gamma (X_\et, \RR^c (n)) \to
    R\Gamma_\BM (G_\RR, X(\CC), \RR (n)) [1]. \]
\end{definition}

Now we state our conjecture about the regulator, which will play an important
role in everything that follows.

\begin{conjecture}
  $\mathbf{B} (X,n)$: given an arithmetic scheme $X$ with smooth
  quasi-projective $X_\CC$ and $n < 0$, the regulator morphism $Reg_{X,n}$
  induces a quasi-isomorphism of complexes of real vector spaces
  \[ Reg_{X,n}^\vee\colon R\Gamma_c (G_\RR, X (\CC), \RR (n)) [-1] \to
    \RHom (R\Gamma (X_\et, \ZZ^c (n)), \RR). \]
\end{conjecture}

\begin{remark}
  If $X/\FF_q$ is a variety over a finite field, then $X (\CC) = \emptyset$,
  so the regulator map is not interesting. Indeed, in our setting, its purpose
  is to take care of the Archimedian places of $X$. In this case
  $\mathbf{B} (X,n)$ implies that $H^i (X_\et, \ZZ^c (n))$ are torsion groups.
  However, by \cite[Proposition~4.2]{Beshenov-Weil-etale-1},
  Conjecture~$\mathbf{L}^c (X_\et, n)$ already implies that
  $H^i (X_\et, \ZZ^c (n))$ are finite groups.
\end{remark}

\begin{remark}
  \label{rmk:regulator-is-defined-for-XC-smooth-quasi-proj}
  We reiterate that our construction of $Reg_{X,n}$ works for $X_{\red,\CC}$
  smooth quasi-projective. In everything that follows, whenever the regulator
  morphism or Conjecture~$\mathbf{B} (X,n)$ is brought, we tacitly assume this
  restriction. This is rather unfortunate, since Weil-\'{e}tale cohomology was
  constructed in \cite{Beshenov-Weil-etale-1} for any arithmetic scheme,
  assuming only Conjecture~$\mathbf{L}^c (X_\et,n)$. Defining the regulator for
  singular $X_{\red,\CC}$ is an interesting project for future work.
\end{remark}

%%%%%%%%%%%%%%%%%%%%%%%%%%%%%%%%%%%%%%%%%%%%%%%%%%%%%%%%%%%%%%%%%%%%%%%%%%%%%%%%

\section{Vanishing order Conjecture $\mathbf{VO} (X,n)$}
\label{sec:vanishing-order-conjecture}

Assuming that $\zeta (X,s)$ admits a meromorphic continuation around
$s = n < 0$, we make the following conjecture for the vanishing order at
$s = n$.

\begin{conjecture}
  $\mathbf{VO} (X,n)$: one has
  \[ \ord_{s=n} \zeta (X,s) =
    \sum_{i \in \ZZ} (-1)^i \cdot i \cdot \rk_\ZZ H^i_\Wc (X, \ZZ (n)). \]
\end{conjecture}

We note that the right-hand side makes sense under
Conjecture~$\mathbf{L}^c (X_\et,n)$, which implies that $H^i_\Wc (X, \ZZ (n))$
are finitely generated groups, trivial for $|i| \gg 0$;
see \cite[Proposition~7.12]{Beshenov-Weil-etale-1}.

\begin{remark}
  Conjecture~$\mathbf{VO} (X,n)$ is similar to
  \cite[Conjecture~5.11]{Flach-Morin-2018}. If $X$ is proper and regular, then
  $\mathbf{VO} (X,n)$ is the same as Flach and Morin's vanishing order
  conjecture. Indeed, the latter is
  \begin{equation}
    \label{eqn:FM-vanishing-order}
    \ord_{s = n} \zeta (X,s) =
    \sum_{i\in \ZZ} (-1)^i \cdot i \cdot \dim_\RR H^i_{\ar,c} (X, \widetilde{\RR}(n)),
  \end{equation}
  where
  \[ R\Gamma_{\ar,c} (X, \widetilde{\RR}(n)) \dfn
    R\Gamma_c (X, \RR(n)) \oplus R\Gamma_c (X, \RR(n)) [-1]. \]
  Moreover, \cite[Proposition 4.14]{Flach-Morin-2018}, gives a distinguished
  triangle
  \begin{multline*}
    R\Gamma_\dR (X_\RR/\RR) / \Fil^n [-2] \to
    R\Gamma_{\ar,c} (X, \widetilde{\RR}(n)) \to
    R\Gamma_\Wc (X, \ZZ(n)) \otimes \RR \\
    \to R\Gamma_\dR (X_\RR/\RR) / \Fil^n [-1]
  \end{multline*}
  So, in case of $n < 0$ we have
  $R\Gamma_{\ar,c} (X, \widetilde{\RR}(n)) \cong
  R\Gamma_\Wc (X, \ZZ(n)) \otimes \RR$ and
  \eqref{eqn:FM-vanishing-order} is exactly Conjecture~$\mathbf{VO} (X,n)$.
\end{remark}

\begin{remark}
  The alternating sum in Conjecture~$\mathbf{VO} (X,n)$ is the so-called
  \textbf{secondary Euler characteristic}
  \[ \chi' (R\Gamma_\Wc (X, \ZZ (n))) \dfn
    \sum_{i \in \ZZ} (-1)^i \cdot i \cdot \rk_\ZZ H^i_\Wc (X, \ZZ (n)). \]
  The calculations below show that the usual Euler characteristic of
  $R\Gamma_\Wc (X, \ZZ (n))$ vanishes, assuming
  Conjectures~$\mathbf{L}^c (X_\et,n)$ and $\mathbf{B} (X,n)$.  See
  \cite{Ramachandran-2016} for more details on the secondary Euler
  characteristic and its occurrences in nature.
\end{remark}

Under the regulator conjecture, our vanishing order formula takes the form of
the usual Euler characteristic of equivariant cohomology
$R\Gamma_c (G_\RR, X(\CC), \RR (n))$ or motivic cohomology
$R\Gamma (X_\et, \ZZ^c (n)) [1]$.

\begin{proposition}
  \label{prop:VO(X,n)-assuming-B(X,n)}
  Assuming $\mathbf{L}^c (X_\et, n)$ and $\mathbf{B} (X,n)$,
  Conjecture~$\mathbf{VO} (X,n)$ is equivalent to
  \begin{align*}
    \ord_{s=n} \zeta (X,s) & = \chi (R\Gamma_c (G_\RR, X(\CC), \RR (n))
    = \sum_{i \in \ZZ} (-1)^i \dim_\RR H^i_c (X(\CC), \RR (n))^{G_\RR} \\
                           & = -\chi (R\Gamma (X_\et, \ZZ^c (n)))
    = \sum_{i \in \ZZ} (-1)^{i+1} \rk_\ZZ H^i (X_\et, \ZZ^c (n)).
  \end{align*}
  Moreover, we have
  $$\chi (R\Gamma_\Wc (X, \ZZ(n))) = 0.$$

  \begin{proof}
    Thanks to \cite[Proposition~7.13]{Beshenov-Weil-etale-1}, the Weil-\'{e}tale
    complex tensored with $\RR$ splits as
    \[ R\Gamma_\Wc (X,\RR (n)) \cong
      \RHom (R\Gamma (X_\et, \ZZ^c (n)), \RR) [-1] \oplus
      R\Gamma_c (G_\RR, X (\CC), \RR (n)) [-1]. \]
    Assuming Conjecture~$\mathbf{B} (X,n)$, we also have a quasi-isomorphism
    \[ R\Gamma_c (G_\RR, X (\CC), \RR (n)) [-1] \cong
      \RHom (R\Gamma (X_\et, \ZZ^c (n)), \RR), \]
    so that
    \[ \dim_\RR H^i_\Wc (X,\RR(n)) =
      \dim_\RR H^{i-1}_c (X (\CC), \RR (n))^{G_\RR} +
      \dim_\RR H^{i-2}_c (X (\CC), \RR (n))^{G_\RR}. \]
    Thus, we can rewrite the sum
    \begin{align*}
      \sum_{i \in \ZZ} (-1)^i \cdot i \cdot \rk_\ZZ H^i_\Wc (X, \ZZ (n)) & = \sum_{i \in \ZZ} (-1)^i \cdot i \cdot \dim_\RR H^i_\Wc (X, \RR (n)) \\
                                                                         & = \sum_{i \in \ZZ} (-1)^i \cdot i \cdot
                                                                           \dim_\RR H^{i-1}_c (X (\CC), \RR (n))^{G_\RR} \\
      & \quad\quad + \sum_{i \in \ZZ} (-1)^i \cdot i \cdot \dim_\RR H^{i-2}_c (X (\CC), \RR (n))^{G_\RR} \\
                                                                         & = -\sum_{i \in \ZZ} (-1)^i \, \dim_\RR H^{i-1}_c (X (\CC), \RR (n))^{G_\RR} \\
                                                                         & = \chi (R\Gamma_c (G_\RR, X (\CC), \RR (n)).
    \end{align*}
    Similarly,
    \begin{align*}
      \sum_{i \in \ZZ} (-1)^i \cdot i \cdot \rk_\ZZ H^i_\Wc (X, \ZZ (n)) & = \chi (\RHom (R\Gamma (X_\et, \ZZ^c (n)), \RR) [1]) \\
                                                                         & = -\chi (R\Gamma (X_\et, \ZZ^c (n))).
    \end{align*}
    These considerations also show that the usual Euler characteristic of
    $R\Gamma_\Wc (X, \ZZ(n))$ vanishes.
  \end{proof}
\end{proposition}

\begin{remark}
  Conjecture~$\mathbf{VO} (X,n)$ is related to a conjecture of Soul\'{e}
  \cite[Conjecture~2.2]{Soule-1984-ICM}, which originally reads in terms of
  $K'$-theory
  \[ \ord_{s=n} \zeta (X,s) =
    \sum_{i \in \ZZ} (-1)^{i+1} \, \dim_\QQ K'_i (X)_{(n)}. \]
  As explained in \cite[Remark~43]{Kahn-2005}, this can be rewritten in
  terms of Borel--Moore motivic homology as
  $$\sum_{i \in \ZZ} (-1)^{i+1} \, \dim_\QQ H_i^{BM} (X, \QQ (n)).$$
  In our setting, $H^i (X_\et, \ZZ^c (n))$ plays the role of Borel--Moore
  homology, which explains the formula
  \[ \ord_{s=n} \zeta (X,s) =
    \sum_{i \in \ZZ} (-1)^{i+1} \rk_\ZZ H^i (X_\et, \ZZ^c (n)). \]
\end{remark}

\begin{remark}[{\cite[Proposition~5.13]{Flach-Morin-2018}}]
  \label{rmk:archimedian-euler-factor}
  As for the formula
  \[ \ord_{s=n} \zeta (X,s) =
    \sum_{i \in \ZZ} (-1)^i \dim_\RR H^i_c (X(\CC), \RR (n))^{G_\RR}, \]
  it essentially means that the vanishing order at $s = n < 0$ comes from the
  Archimedian $\Gamma$-factor appearing in the (hypothetical) functional
  equation, as explained in \cite[\S\S 3,4]{Serre-1970}
  (see also \cite[\S 4]{Flach-Morin-2020}).

  Indeed, under the assumption that $X_\CC$ is a smooth projective variety, we
  consider the Hodge decomposition
  \[ H^i (X (\CC), \CC) = \bigoplus_{p+q = i} H^{p,q}, \]
  which carries an action of $G_\RR = \{ id, \sigma \}$ such that
  $\sigma (H^{p,q}) = H^{q,p}$. We set $h^{p,q} = \dim_\CC H^{p,q}$.
  For $p = i/2$ we consider the eigenspace decomposition
  $H^{p,p} = H^{p,+} \oplus H^{p,-}$, where
  \begin{align*}
    H^{p,+} & = \{ x \in H^{p,p} \mid \sigma (x) = (-1)^p\,x \},\\
    H^{p,-} & = \{ x \in H^{p,p} \mid \sigma (x) = (-1)^{p+1}\,x \},
  \end{align*}
  and set $h^{p,\pm} = \dim_\CC H^{p,\pm}$ accordingly.
  The completed zeta function
  $$\zeta (\overline{X}, s) = \zeta (X, s)\,\zeta (X_\infty, s)$$
  is expected to satisfy a functional equation of the form
  \[ A^{\frac{d-s}{2}}\,\zeta (\overline{X},d-s) =
    A^{\frac{s}{2}}\,\zeta (\overline{X},s). \]
  Here
  \begin{gather*}
    \zeta (X_\infty, s) = \prod_{i\in \ZZ} L_\infty (H^i (X),s)^{(-1)^i}, \\
    L_\infty (H^i (X), s) =
    \prod_{p = i/2} \Gamma_\RR (s - p)^{h^{p,+}}\,\Gamma_\RR (s-p+1)^{h^{p,-}} \,
    \prod_{\substack{p + q = i \\ p < q}} \Gamma_\CC (s - p)^{h^{p,q}}, \\
    \Gamma_\RR (s) = \pi^{-s/2} \, \Gamma (s/2), \quad
    \Gamma_\CC (s) = (2\pi)^{-s} \, \Gamma (s).
  \end{gather*}

  Therefore, the expected vanishing order at $s = n < 0$ is
  \begin{align*}
    \ord_{s=n} \zeta (X,s) & = -\ord_{s=n} \zeta (X_\infty,s) \\
                           & = -\sum_{i\in \ZZ} (-1)^i \ord_{s=n} L_\infty (H^i (X), s) \\
                           & = \sum_{i\in \ZZ} (-1)^i \Bigl(\sum_{p = i/2} h^{p,(-1)^{n-p}} +
    \sum_{\substack{p + q = i \\ p < q}} h^{p,q}\Bigr).
  \end{align*}
  The last equality follows from the fact that $\Gamma (s)$ has simple poles at
  all $s = n \le 0$. We have
  \begin{align*}
    \dim_\RR H^i (X (\CC), \RR (n))^{G_\RR} & = \dim_\RR H^i (X (\CC), \RR)^{\sigma = (-1)^n} \\
                                            & = \dim_\CC H^i (X (\CC), \CC)^{\sigma = (-1)^n} \\
                                            & = \sum_{p = i/2} h^{p,(-1)^{n-p}} + \sum_{\substack{p + q = i \\ p < q}} h^{p,q}.
  \end{align*}
  Here the terms $h^{p,q}$ with $p < q$ come from $\sigma (H^{p,q}) = H^{q,p}$,
  while $h^{p,(-1)^{n-p}}$ come from the action on $H^{p,p}$.
  We see that our conjectural formula recovers the expected vanishing order.
\end{remark}

Let us look at some particular examples when the meromorphic continuation for
$\zeta (X,s)$ is known.

\begin{example}
  \label{example:VO(X,n)-for-number-rings}
  Suppose that $X = \Spec \mathcal{O}_F$ is the spectrum of the ring of integers
  of a number field $F/\QQ$. Let $r_1$ be the number of real embeddings
  $F \hookrightarrow \RR$ and $r_2$ be the number of conjugate pairs of complex
  embeddings $F \hookrightarrow \CC$. The space $X (\CC)$ with the action of
  complex conjugation can be visualized as follows:
  
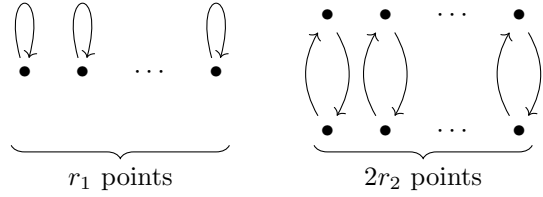
\[ \begin{tikzpicture}
    \matrix(m)[matrix of math nodes, row sep=1em, column sep=1em,
    text height=1ex, text depth=0.2ex]{
      ~ & ~ & ~ & ~ & ~ & \bullet & \bullet & \cdots & \bullet \\
      \bullet & \bullet & \cdots & \bullet \\
      ~ & ~ & ~ & ~ & ~ & \bullet & \bullet & \cdots & \bullet \\};

    \draw[->] (m-2-1) edge[loop above,min distance=10mm] (m-2-1);
    \draw[->] (m-2-2) edge[loop above,min distance=10mm] (m-2-2);
    \draw[->] (m-2-4) edge[loop above,min distance=10mm] (m-2-4);

    \draw[->] (m-1-6) edge[bend left] (m-3-6);
    \draw[->] (m-1-7) edge[bend left] (m-3-7);
    \draw[->] (m-1-9) edge[bend left] (m-3-9);

    \draw[->] (m-3-6) edge[bend left] (m-1-6);
    \draw[->] (m-3-7) edge[bend left] (m-1-7);
    \draw[->] (m-3-9) edge[bend left] (m-1-9);

    \draw [decorate,decoration={brace,amplitude=5pt,mirror}] ($(m-3-1)+(-0.5em,-0.5em)$) -- ($(m-3-4)+(0.5em,-0.5em)$);
    \draw [decorate,decoration={brace,amplitude=5pt,mirror}] ($(m-3-6)+(-0.5em,-0.5em)$) -- ($(m-3-9)+(0.5em,-0.5em)$);

    \draw ($(m-3-1)!.5!(m-3-4)$) node[yshift=-2em,anchor=base] {$r_1$ points};
    \draw ($(m-3-6)!.5!(m-3-9)$) node[yshift=-2em,anchor=base] {$2 r_2$ points};
  \end{tikzpicture} \]

  The complex $R\Gamma_c (X (\CC), \RR (n))$ consists of a single $G_\RR$-module
  in degree $0$ given by
  $$\RR (n)^{\oplus r_1} \oplus (\RR (n) \oplus \RR (n))^{\oplus r_2},$$
  with the action of $G_\RR$ on the first summand $\RR (n)^{\oplus r_1}$ via the
  complex conjugation and the action on the second summand
  $(\RR (n) \oplus \RR (n))^{\oplus r_2}$ via
  $(x,y) \mapsto (\overline{y}, \overline{x})$. The corresponding real space of
  fixed points has dimension
  \[ \dim_\RR H^0_c (G_\RR, X (\CC), \RR (n)) = \begin{cases}
      r_2, & n \text{ odd},\\
      r_1 + r_2, & n \text{ even},\\
    \end{cases} \]
  which indeed coincides with the vanishing order of the Dedekind zeta function
  $\zeta (X,s) = \zeta_F (s)$ at $s = n < 0$.

  On the motivic cohomology side, for $n < 0$ the groups
  $H^i (X_\et, \ZZ^c (n))$ are finite, except for $i = -1$, where by
  \cite[Proposition~4.14]{Geisser-2017}
  \[ \rk_\ZZ H^{-1} (X_\et, \ZZ^c (n)) = \begin{cases}
      r_2, & n \text{ odd},\\
      r_1 + r_2, & n \text{ even}.
    \end{cases} \]
\end{example}

\begin{example}
  Suppose that $X$ is a variety over a finite field $\FF_q$. Then the vanishing
  order conjecture is not very interesting, because the formula yields
  \begin{align*}
    \ord_{s=n} \zeta (X,s) & = \sum_{i \in \ZZ} (-1)^i \dim_\RR H^i_c (X(\CC), \RR (n))^{G_\RR} \\
                           & = \sum_{i \in \ZZ} (-1)^{i+1} \rk_\ZZ H^i (X_\et, \ZZ^c (n)) = 0,
  \end{align*}
  since $X (\CC) = \emptyset$, and also because $\mathbf{L}^c (X_\et, n)$
  implies $\rk_\ZZ H^i (X_\et, \ZZ^c (n)) = 0$ for all $i$ in the case of
  varieties over finite fields, as observed in
  \cite[Proposition~4.2]{Beshenov-Weil-etale-1}. Therefore, the conjecture
  simply asserts that $\zeta (X,s)$ has no zeros or poles at $s = n < 0$.  This
  is indeed the case. We have $\zeta (X,s) = Z (X,q^{-s})$, where
  $$Z (X,t) = \exp \Bigl(\sum_{k\ge 1} \frac{\# X (\FF_{q^k})}{k}\,t^k\Bigr)$$
  is the Hasse--Weil zeta function. According to Deligne's work on Weil's
  conjectures \cite{Deligne-Weil-II}, the zeros and poles of $Z (X,s)$ satisfy
  $|s| = q^{-w/2}$, where $0 \le w \le 2 \dim X$
  (see e.g. \cite[pp.\,26--27]{Katz-1994}). In particular, $q^{-s}$ for
  $s = n < 0$ is neither a zero nor a pole of $Z (X,s)$.

  We also note that our definition of $H^i_\Wc (X, \ZZ(n))$, and pretty much
  everything said above, only makes sense for $n < 0$. Already for $n = 0$, for
  example, the zeta function of a smooth projective curve $X/\FF_q$ has a simple
  pole at $s = 0$.
\end{example}

\begin{example}
  Let $X = E$ be an integral model of an elliptic curve over $\QQ$. Then, as a
  consequence of the modularity theorem
  (Wiles--Breuil--Conrad--Diamond--Taylor), it is known that $\zeta (E,s)$
  admits a meromorphic continuation satisfying the functional equation with the
  $\Gamma$-factors discussed in Remark~\ref{rmk:archimedian-euler-factor}.
  In this particular case $\ord_{s=n} \zeta (E,s) = 0$ for all $n < 0$. This is
  consistent with the fact that
  $\chi (R\Gamma_c (G_\RR, E (\CC), \RR (n))) = 0$.

  Indeed, the equivariant cohomology groups $H^i_c (E (\CC), \RR (n))^{G_\RR}$
  are the following:
  \begin{center}
    \renewcommand{\arraystretch}{1.5}
    \begin{tabular}{rccc}
      \hline
      & $i = 0$ & $i = 1$ & $i = 2$ \\
      \hline
      $n$ even: & $\RR$ & $\RR$ & $0$ \\
      $n$ odd: & $0$ & $\RR$ & $\RR$ \\
      \hline
    \end{tabular}
  \end{center}
  ---see, for example, the calculation in \cite[Lemma~A.6]{Siebel-2019}.
\end{example}

%%%%%%%%%%%%%%%%%%%%%%%%%%%%%%%%%%%%%%%%%%%%%%%%%%%%%%%%%%%%%%%%%%%%%%%%%%%%%%%%

\section{Special value Conjecture~$\mathbf{C} (X,n)$}
\label{sec:special-value-conjecture}

\begin{definition}
  We define a morphism of complexes
  \[ \smile\theta\colon R\Gamma_\Wc (X,\ZZ(n)) \otimes \RR \to
    R\Gamma_\Wc (X,\ZZ(n)) [1] \otimes \RR \]
  using the splitting \cite[Proposition~7.13]{Beshenov-Weil-etale-1}
  \[ R\Gamma_\Wc (X, \RR (n)) \cong
    \RHom (R\Gamma (X_\et, \ZZ^c (n)), \RR) [-1] \oplus
    R\Gamma_c (G_\RR, X (\CC), \RR (n)) [-1] \]
  as follows:
  \[ \begin{tikzcd}
      R\Gamma_\Wc (X, \RR(n)) \ar{d}{\cong}\ar[dashed]{r}{\smile\theta} & R\Gamma_\Wc (X, \RR(n)) [1]\ar{d}{\cong} \\
      \RHom (R\Gamma (X_\et, \ZZ^c (n)), \RR) [-1] & \RHom (R\Gamma (X_\et, \ZZ^c (n)), \RR) \\[-2em]
      \oplus & \oplus \\[-2em]
      R\Gamma_c (G_\RR, X (\CC), \RR (n)) [-1]\ar{uur}[description]{Reg_{X,n}^\vee} & R\Gamma_c (G_\RR, X (\CC), \RR (n))
    \end{tikzcd} \]
\end{definition}

\begin{lemma}
  \label{lemma:smile-theta}
  Assuming Conjectures $\mathbf{L}^c (X_\et,n)$ and $\mathbf{B} (X,n)$, the
  morphism $\smile\theta$ induces a long exact sequence of finite dimensional
  real vector spaces
  \[ \cdots \to H^{i-1}_\Wc (X, \RR (n))
    \xrightarrow{\smile\theta}
    H^i_\Wc (X, \RR (n))
    \xrightarrow{\smile\theta}
    H^{i+1}_\Wc (X, \RR (n)) \to \cdots \]

  \begin{proof}
    We obtain a sequence
    \[ \begin{tikzcd}[column sep=1em]
        \cdots\ar{r} & H^i_\Wc (X, \RR(n))\ar{d}{\cong}\ar[dashed]{r}{\smile\theta} & H^{i+1}_\Wc (X, \RR(n))\ar{d}{\cong} \ar{r} & \cdots \\
        & \Hom (H^{-i-1} (X_\et, \ZZ^c (n)), \RR) & \Hom (H^{-i-2} (X_\et, \ZZ^c (n)), \RR) \\[-2em]
        \cdots & \oplus & \oplus & \cdots \\[-2em]
        & H_c^{i-1} (G_\RR, X (\CC), \RR (n))\ar{uur}[description]{\cong} & H_c^i (G_\RR, X (\CC), \RR (n))
      \end{tikzcd} \]
    The diagonal arrows are isomorphisms according to $\mathbf{B} (X,n)$,
    so the sequence is exact.
  \end{proof}
\end{lemma}

The Weil-\'{e}tale complex $R\Gamma_\Wc (X, \ZZ(n))$ is defined in
\cite[\S 7]{Beshenov-Weil-etale-1} up to a \emph{non-unique} isomorphism in the
derived category $\mathbf{D} (\ZZ)$ via a distinguished triangle
\begin{equation}
  \label{eqn:triangle-defining-RGamma-Wc}
  R\Gamma_\Wc (X, \ZZ(n)) \to R\Gamma_\fg (X,\ZZ(n)) \xrightarrow{i_\infty}
  R\Gamma_c (G_\RR, X (\CC), \ZZ (n)) \to [1]
\end{equation}
This is rather awkward, and there should be a better, more canonical
construction of $R\Gamma_\Wc (X, \ZZ(n))$. For our purposes, however, this is
not much of a problem, since the special value conjecture is not formulated in
terms of $R\Gamma_\Wc (X, \ZZ (n))$, but in terms of its determinant
$\det_\ZZ R\Gamma_\Wc (X, \ZZ(n))$ (see Appendix~\ref{app:determinants}), which
is well-defined.

\begin{lemma}
  \label{lemma:determinant-of-RGamma-Wc-well-defined}
  The determinant $\det_\ZZ R\Gamma_\Wc (X, \ZZ(n))$ is defined up to a
  canonical isomorphism.

  \begin{proof}
    Two different choices for the mapping fiber in
    \eqref{eqn:triangle-defining-RGamma-Wc} yield an isomorphism of
    distinguished triangles
    \[ \begin{tikzcd}
        R\Gamma_\Wc (X, \ZZ(n)) \ar{r}\ar{d}{f}[swap]{\cong} & R\Gamma_\fg (X,\ZZ(n)) \ar{r}{i_\infty}\ar{d}{id} & R\Gamma_c (G_\RR, X (\CC), \ZZ (n)) \ar{r}\ar{d}{id} & {[1]}\ar{d}{f}[swap]{\cong} \\
        R\Gamma_\Wc (X, \ZZ(n))' \ar{r} & R\Gamma_\fg (X,\ZZ(n)) \ar{r}{i_\infty} & R\Gamma_c (G_\RR, X (\CC), \ZZ (n)) \ar{r} & {[1]}
      \end{tikzcd} \]
    The idea is to use functoriality of determinants with respect to
    isomorphisms of distinguished triangles
    (see Appendix~\ref{app:determinants}). The only technical problem is that
    whenever $X (\RR) \ne \emptyset$, the complexes $R\Gamma_\fg (X,\ZZ(n))$ and
    $R\Gamma_c (G_\RR, X (\CC), \ZZ (n))$ are not perfect, but may have finite
    $2$-torsion in $H^i (-)$ for arbitrarily big $i$ (in
    \cite{Beshenov-Weil-etale-1} we called such complexes \textbf{almost
      perfect}). On the other hand, the determinants are defined only for
    perfect complexes. Fortunately, $H^i (i_\infty^*)$ is an isomorphism for
    $i \gg 0$ by the boundedness of $H^i_\Wc (X, \ZZ(n))$
    \cite[Proposition~7.12]{Beshenov-Weil-etale-1}, so that for $m$ big enough
    we can take the corresponding canonical truncations $\tau_{\le m}$:
    \[ \begin{tikzcd}[row sep=1.5em, column sep=0.75em, font=\small]
        \tau_{\le m} R\Gamma_\Wc (X, \ZZ(n)) \ar{r}\ar{d}{\cong} & \tau_{\le m} R\Gamma_\fg (X,\ZZ(n)) \ar{r}\ar{d} & \tau_{\le m} R\Gamma_c (G_\RR, X (\CC), \ZZ (n)) \ar{r}\ar{d} & {[1]}\ar{d}{\cong} \\
        R\Gamma_\Wc (X, \ZZ(n)) \ar{r}\ar{d} & R\Gamma_\fg (X,\ZZ(n)) \ar{r}{i_\infty}\ar{d} & R\Gamma_c (G_\RR, X (\CC), \ZZ (n)) \ar{r}\ar{d} & {[1]}\ar{d} \\
        0 \ar{r}\ar{d} & \tau_{\ge m+1} R\Gamma_\fg (X,\ZZ(n)) \ar{r}{\cong}\ar{d} & \tau_{\ge m+1} R\Gamma_c (G_\RR, X (\CC), \ZZ (n)) \ar{r}\ar{d} & 0\ar{d} \\
        {[1]} \ar{r} & {[1]} \ar{r} & {[1]} \ar{r} & {[2]}
      \end{tikzcd} \]

    The truncations give us (rotating the triangles)
    \[ \begin{tikzcd}[column sep=1em,font=\small]
        \tau_{\le m} R\Gamma_c (G_\RR, X (\CC), \ZZ (n))[-1] \ar{r}\ar{d}{id} & R\Gamma_\Wc (X, \ZZ (n))\ar{d}{f}[swap]{\cong}\ar{r} & \tau_{\le m} R\Gamma_\fg (X,\ZZ(n)) \ar{d}{id}\ar{r} & {[0]} \ar{d}{id} \\
        \tau_{\le m} R\Gamma_c (G_\RR, X (\CC), \ZZ (n))[-1] \ar{r} & R\Gamma_\Wc (X, \ZZ (n))' \ar{r} & \tau_{\le m} R\Gamma_\fg (X,\ZZ(n)) \ar{r} & {[0]}
      \end{tikzcd} \]
    By the functoriality of determinants with respect to isomorphisms of
    distinguished triangles (see Appendix~\ref{app:determinants}), we have a
    commutative diagram
    \[ \begin{tikzcd}
        \begin{array}{c} \det_\ZZ \tau_{\le m} R\Gamma_c (G_\RR, X (\CC), \ZZ (n))[-1] \\ \otimes \\ \det_\ZZ \tau_{\le m} R\Gamma_\fg (X,\ZZ(n)) \end{array} \ar{r}{i}[swap]{\cong} \ar{d}{id} & \det_\ZZ R\Gamma_\Wc (X,\ZZ(n)) \ar{d}{\det_\ZZ (f)}[swap]{\cong} \\
        \begin{array}{c} \det_\ZZ \tau_{\le m} R\Gamma_c (G_\RR, X (\CC), \ZZ (n))[-1] \\ \otimes \\ \det_\ZZ \tau_{\le m} R\Gamma_\fg (X, \ZZ (n)) \end{array} \ar{r}{i'}[swap]{\cong} & \det_\ZZ R\Gamma_\Wc (X,\ZZ(n))
      \end{tikzcd} \]
    so that $\det_\ZZ (f) = i'\circ i^{-1}$.
  \end{proof}
\end{lemma}

\begin{lemma}
  The non-canonical splitting
  \[ R\Gamma_\Wc (X, \RR(n)) \cong
    \RHom (R\Gamma (X_\et, \ZZ^c (n)), \RR) [-1] \oplus
    R\Gamma_c (G_\RR, X (\CC), \RR (n)) [-1] \]
  from \cite[Proposition~7.13]{Beshenov-Weil-etale-1} yields a canonical
  isomorphism of determinants
  \[ \det_\RR R\Gamma_\Wc (X, \RR (n)) \cong
    \begin{array}{c}
      \det_\RR \RHom (R\Gamma (X_\et, \ZZ^c (n)), \RR) [-1] \\
      \otimes_\RR \\
      \det_\RR R\Gamma_c (G_\RR, X (\CC), \RR (n)) [-1]
    \end{array} \]

  \begin{proof}
    This is similar to the previous lemma; in fact, after tensoring with $\RR$,
    we obtain perfect complexes of real vector spaces, so the truncations are no
    longer needed. By \cite[Proposition~7.4]{Beshenov-Weil-etale-1} we have
    $i_\infty^* \otimes \RR = 0$, so there is an isomorphism of triangles
    \begin{equation}
      \label{eqn:splitting-of-RGamma-Wc-triangles}
      \begin{tikzcd}[row sep=1em]
        R\Gamma_c (G_\RR, X (\CC), \RR (n)) [-1] \ar{r}{id}\ar{d} & R\Gamma_c (G_\RR, X (\CC), \RR (n)) [-1]\ar{d} \\
        R\Gamma_\Wc (X, \RR (n)) \ar[dashed]{r}{f}[swap]{\cong}\ar{d} & \begin{array}{c} \RHom (R\Gamma (X_\et, \ZZ^c (n)), \RR) [-1] \\ \oplus \\ R\Gamma_c (G_\RR, X (\CC), \RR (n)) [-1] \end{array}\ar{d} \\
        R\Gamma_\fg (X, \RR (n)) \ar{r}{g\otimes \RR}[swap]{\cong}\ar{d} & \RHom (R\Gamma (X_\et, \ZZ^c (n)), \RR) [-1] \ar{d} \\
        R\Gamma_c (G_\RR, X (\CC), \RR (n)) \ar{r}{id} & R\Gamma_c (G_\RR, X (\CC), \RR (n))
      \end{tikzcd}
    \end{equation}
    Here the third horizontal arrow comes from the triangle defining
    $R\Gamma_\fg (X, \ZZ(n))$:
    \begin{multline*}
      \RHom (R\Gamma (X_\et, \ZZ^c (n)), \RR) [-2] \xrightarrow{\alpha_{X,n}}
      R\Gamma_c (X_\et, \ZZ (n)) \to
      R\Gamma_\fg (X, \ZZ(n)) \\
      \xrightarrow{g} \RHom (R\Gamma (X_\et, \ZZ^c (n)), \RR) [-1]
    \end{multline*}
    tensored with $\RR$ (see \cite[Proposition~5.7]{Beshenov-Weil-etale-1}).
    The distinguished column on the right-hand side of
    \eqref{eqn:splitting-of-RGamma-Wc-triangles} is the direct sum
    \[ \begin{tikzcd}[row sep=1.5em]
        R\Gamma_c (G_\RR, X (\CC), \RR (n)) [-1]\ar{d}{id} &[-2em] &[-2em] 0\ar{d} \\
        R\Gamma_c (G_\RR, X (\CC), \RR (n)) [-1]\ar{d} & \oplus & \RHom (R\Gamma (X_\et, \ZZ^c (n)), \RR) [-1]\ar{d}{id} \\
        0 \ar{d} & & \RHom (R\Gamma (X_\et, \ZZ^c (n)), \RR) [-1] \ar{d} \\
        R\Gamma_c (G_\RR, X (\CC), \RR (n)) & & 0
      \end{tikzcd} \]
    The splitting isomorphism $f$ in
    \eqref{eqn:splitting-of-RGamma-Wc-triangles} is not canonical at
    all. However, after taking the determinants, we obtain a commutative diagram
    (see Appendix~\ref{app:determinants})
    \[ \begin{tikzcd}[column sep=1em, font=\small]
        \begin{array}{c} \det_\RR R\Gamma_c (G_\RR, X (\CC), \RR (n)) [-1] \\ \otimes_\RR \\ \det_\RR R\Gamma_\fg (X, \RR(n)) \end{array} \ar{r}{i}[swap]{\cong} \ar{d}{id \otimes \det_\RR (g\otimes \RR)}[swap]{\cong} & \det_\RR R\Gamma_\Wc (X, \RR (n)) \ar{d}{\det_\RR (f)}[swap]{\cong}\ar[dashed]{dl}{\cong} \\
        \begin{array}{c} \det_\RR R\Gamma_c (G_\RR, X (\CC), \RR (n)) [-1] \\ \otimes_\RR \\ \det_\RR \RHom (R\Gamma (X_\et, \ZZ^c (n)), \RR) [-1] \end{array} \ar{r}{i'}[swap]{\cong} & \det_\RR \left(\!\!\!\begin{array}{c} \RHom (R\Gamma (X_\et, \ZZ^c (n)), \RR) [-1] \\ \oplus \\ R\Gamma_c (G_\RR, X (\CC), \RR (n)) [-1] \end{array}\!\!\!\right)
      \end{tikzcd} \]
    The dashed diagonal arrow is the desired canonical isomorphism.
  \end{proof}
\end{lemma}

\begin{definition}
  Given an arithmetic scheme $X$ and $n < 0$, assume
  Conjectures~$\mathbf{L}^c (X_\et, n)$ and $\mathbf{B} (X,n)$. Consider the
  quasi-isomorphism
  \begin{multline}\small
    \label{eqn:definition-of-lambda}
    \left(\!\!\!\begin{array}{c} R\Gamma_c (G_\RR, X (\CC), \RR (n)) [-2] \\ \oplus \\ R\Gamma_c (G_\RR, X (\CC), \RR (n)) [-1] \end{array}\!\!\!\right)
    \xrightarrow[\cong]{Reg_{X,n}^\vee [-1] \oplus id}
    \left(\!\!\!\begin{array}{c} \RHom (R\Gamma (X_\et, \ZZ^c (n)), \RR) [-1] \\ \oplus \\ R\Gamma_c (G_\RR, X (\CC), \RR (n)) [-1] \end{array}\!\!\!\right) \\
    \\
    \xrightarrow[\cong]{\text{split}} R\Gamma_\Wc (X, \RR (n))
  \end{multline}
  Note that the first complex has determinant
  \[ \det_\RR \left(\!\!\begin{array}{c} R\Gamma_c (G_\RR, X (\CC), \RR (n)) [-2] \\ \oplus \\ R\Gamma_c (G_\RR, X (\CC), \RR (n)) [-1] \end{array}\!\!\right) \cong
    \begin{array}{c} \det_\RR R\Gamma_c (G_\RR, X (\CC), \RR (n)) \\ \otimes_\RR \\ (\det_\RR R\Gamma_c (G_\RR, X (\CC), \RR (n)))^{-1} \end{array} \cong \RR, \]
  and for the last complex in \eqref{eqn:definition-of-lambda}, by the
  compatibility with base change, we have a canonical isomorphism
  \[ \det_\RR R\Gamma_\Wc (X, \RR (n)) \cong
    (\det_\ZZ R\Gamma_\Wc (X, \ZZ (n))) \otimes \RR. \]
  Therefore, after taking the determinants, the quasi-isomorphism
  \eqref{eqn:definition-of-lambda} induces a canonical isomorphism
  \begin{equation}
    \label{eqn:morphism-lambda}
    \lambda = \lambda_{X,n}\colon \RR \xrightarrow{\cong}
    (\det_\ZZ R\Gamma_\Wc (X, \ZZ (n))) \otimes \RR.
  \end{equation}
\end{definition}

\begin{remark}
  An equivalent way to define $\lambda$ is
  \begin{multline*}
    \lambda\colon \RR \xrightarrow{\cong}
    \bigotimes_{i\in \ZZ} (\det_\RR H^i_\Wc (X, \RR (n)))^{(-1)^i} \\
    \xrightarrow{\cong} \Bigl(\bigotimes_{i\in \ZZ} (\det_\ZZ H^i_\Wc (X, \ZZ (n)))^{(-1)^i}\Bigr) \otimes \RR \\
    \xrightarrow{\cong} (\det_\ZZ R\Gamma_\Wc (X, \ZZ (n))) \otimes \RR,
  \end{multline*}
  where the first isomorphism comes from Lemma~\ref{lemma:smile-theta}.
\end{remark}

We are ready to state the main conjecture of this paper. The determinant
$\det_\ZZ R\Gamma_\Wc (X, \ZZ (n)))$ is a free $\ZZ$-module of rank $1$, and the
isomorphism \eqref{eqn:morphism-lambda} canonically embeds it in $\RR$. We
conjecture that this embedding gives the special value of $\zeta (X,s)$ at
$s = n$ in the following sense.

\begin{conjecture}
  $\mathbf{C} (X,n)$: let $X$ be an arithmetic scheme and $n < 0$ a strictly
  negative integer. Assuming Conjectures~$\mathbf{L}^c (X_\et, n)$,
  $\mathbf{B} (X,n)$ and the meromorphic continuation of $\zeta (X,s)$ around
  $s = n < 0$, the corresponding special value is determined up to sign by
  \[ \lambda (\zeta^* (X,n)^{-1}) \cdot \ZZ =
    \det_\ZZ R\Gamma_\Wc (X, \ZZ (n)), \]
  where $\lambda$ is the canonical isomorphism \eqref{eqn:morphism-lambda}.
\end{conjecture}

\begin{remark}
  This conjecture is similar to \cite[Conjecture~5.12]{Flach-Morin-2018}.
  When $X$ is proper and regular, the above conjecture is the same as the
  special value conjecture of Flach and Morin, which for $n \in \ZZ$ reads
  \begin{equation}
    \label{eqn:FM-special-value}
    \lambda_\infty \Bigl(\zeta^* (X,n)^{-1} \cdot C (X,n) \cdot \ZZ\Bigr) =
    \Delta (X/\ZZ, n).
  \end{equation}
  Here the fundamental line $\Delta (X/\ZZ,n)$ is defined via
  \[ \Delta (X/\ZZ,n) \dfn
    \det_\ZZ R\Gamma_\Wc (X, \ZZ(n)) \otimes
    \det_\ZZ R\Gamma_\dR (X/\ZZ)/\Fil^n. \]
  If $n < 0$, then
  $\Delta (X/\ZZ,n) = \det_\ZZ R\Gamma_\Wc (X, \ZZ(n))$. Moreover, $C (X,n)$
  in \eqref{eqn:FM-special-value} is a rational number, defined via
  $\prod_p |c_p (X,n)|_p$.  Here $c_p (X,n) \in \QQ_p^\times/\ZZ_p^\times$ are
  the local factors described in \cite[\S 5.4]{Flach-Morin-2018}, but
  \cite[Proposition~5.8]{Flach-Morin-2018} states that if $n \le 0$, then
  $c_p (X,n) \equiv 1 \pmod{\ZZ_p^\times}$ for all $p$. Therefore, $C (X,n) = 1$
  in our situation. Finally, the trivialization isomorphism $\lambda_\infty$ is
  defined exactly as our $\lambda$. Therefore, \eqref{eqn:FM-special-value} for
  $n < 0$ agrees with Conjecture~$\mathbf{C} (X,n)$.

  Flach and Morin prove that their conjecture is consistent with the Tamagawa
  number conjecture of Bloch--Kato--Fontaine--Perrin-Riou
  \cite{Fontaine-Perrin-Riou-1994}; see \cite[\S 5.6]{Flach-Morin-2018} for the
  details.
\end{remark}

\begin{remark}
  Some canonical isomorphisms of determinants involve multiplication by $\pm 1$,
  so it is no surprise that the resulting conjecture is stated up to sign
  $\pm 1$. This is not a major problem, however, since the sign can be recovered
  from the (conjectural) functional equation.
\end{remark}

%%%%%%%%%%%%%%%%%%%%%%%%%%%%%%%%%%%%%%%%%%%%%%%%%%%%%%%%%%%%%%%%%%%%%%%%%%%%%%%%

\section{Case of varieties over finite fields}
\label{sec:finite-fields}

For varieties over finite fields, our special value conjecture corresponds to
the conjectures studied by Geisser in
\cite{Geisser-2004,Geisser-2006,Geisser-2010-arithmetic-homology}.

\begin{proposition}
  \label{prop:C(X,n)-over-finite-fields}
  If $X/\FF_q$ is a variety over a finite field, then under the assumption
  $\mathbf{L}^c (X_\et,n)$, the special value conjecture $\mathbf{C} (X,n)$ is
  equivalent to
  \begin{align}
    \label{eqn:special-value-for-X/Fq}
    \notag \zeta^* (X,n) & = \pm \prod_{i \in \ZZ} |H_\Wc^i (X, \ZZ(n))|^{(-1)^i} \\
                         & = \pm \prod_{i \in \ZZ} |H^i (X_\et, \ZZ^c (n))|^{(-1)^i} \\
    \notag & = \pm \prod_{i \in \ZZ} |H_i^c (X_\ar, \ZZ (n))|^{(-1)^{i+1}},
  \end{align}
  where $H_i^c (X_\ar, \ZZ (n))$ are Geisser's arithmetic homology groups
  defined in \cite{Geisser-2010-arithmetic-homology}.

  \begin{proof}
    Assuming $\mathbf{L}^c (X_\et,n)$, we have, thanks to
    \cite[Proposition~7.7]{Beshenov-Weil-etale-1}
    \[ H^i_\Wc (X, \ZZ (n)) \cong
      \Hom (H^{2-i} (X_\et, \ZZ^c (n)), \QQ/\ZZ) \cong
      \Hom (H_{i-1}^c (X_\ar, \ZZ (n)), \QQ/\ZZ). \]
    The cohomology groups involved are finite and vanish for $|i| \gg 0$ by
    \cite[Proposition~4.2]{Beshenov-Weil-etale-1}, and by
    Lemma~\ref{lemma:determinant-for-torsion-cohomology} the determinant is given by
    \[ \begin{tikzcd}[row sep=0.75em, column sep=0pt]
        \det_\ZZ R\Gamma_\Wc (X, \ZZ (n)) \ar[equals]{d} & \subset & \det_\ZZ R\Gamma_\Wc (X, \ZZ (n)) \otimes \QQ \ar[equals]{d} \\
        \frac{1}{m}\ZZ & \subset & \QQ
      \end{tikzcd} \]
    where
    \[ m = \prod_{i \in \ZZ} |H_\Wc^i (X, \ZZ(n))|^{(-1)^i}. \qedhere \]
  \end{proof}
\end{proposition}

\begin{remark}
  Formulas like \eqref{eqn:special-value-for-X/Fq} were proposed by Lichtenbaum
  early on in \cite{Lichtenbaum-1984}.
\end{remark}

\begin{theorem}
  \label{thm:C(X,n)-over-finite-fields}
  Let $X/\FF_q$ be a variety over a finite field satisfying
  Conjecture~$\mathbf{L}^c (X_\et, n)$ for $n < 0$. Then
  Conjecture~$\mathbf{C} (X,n)$ holds.
\end{theorem}

We note that \eqref{prop:C(X,n)-over-finite-fields} is equivalent to the special
value formula that appears in
\cite[Theorem~4.5]{Geisser-2010-arithmetic-homology}.
Conjecture~$\mathbf{P}_0 (X)$ in the statement of
\cite[Theorem~4.5]{Geisser-2010-arithmetic-homology} is implied by our
Conjecture~$\mathbf{L}^c (X_\et,n)$ thanks to
\cite[Proposition~4.1]{Geisser-2010-arithmetic-homology}. Geisser's proof
eventually reduces to Milne's work \cite{Milne-1986}, but for our case of
$s = n < 0$, the situation is simpler, and we can give a direct explanation,
using earlier results of B\'{a}yer and Neukirch \cite{Bayer-Neukirch-1978}
concerning Grothendieck's trace formula.

\begin{proof}
  By the previous proposition, the conjecture reduces to
  $$\zeta (X,n) = \prod_{i \in \ZZ} |H^i (X_\et, \ZZ^c (n))|^{(-1)^i}.$$

  By duality \cite[Theorem~I]{Beshenov-Weil-etale-1}
  $$|H^{2-i} (X_\et, \ZZ^c (n))| = |H^i_c (X_\et, \ZZ (n))|,$$
  where
  \[ \ZZ (n) \dfn
    \bigoplus_{\ell \ne p} \QQ_\ell/\ZZ_\ell (n) [-1] \dfn
    \bigoplus_{\ell \ne p} \mu_{\ell^\infty}^{\otimes n} [-1] \dfn
    \bigoplus_{\ell \ne p} \varinjlim_r \mu_{\ell^r}^{\otimes n} [-1], \]
  and $p$ is the characteristic of the base field.
  Now $H^i_c (X_\et, \QQ_\ell (n)) = 0$ for $n < 0$, and therefore
  $H^i_c (X_\et, \ZZ_\ell (n)) \cong H^{i-1}_c (X_\et, \QQ_\ell/\ZZ_\ell (n))$.
  This means that our formula can be written as
  \begin{equation}
    \label{eqn:zeta-X/Fq-product-formula}
    \zeta (X,n) =
    \prod_{\ell \ne p} \prod_{i \in \ZZ} |H^i_c (X_\et, \ZZ_\ell (n))|^{(-1)^i}.
  \end{equation}

  Grothendieck's trace formula (see \cite{Grothendieck-FL} or
  \cite[Rapport]{SGA4-1-2}) reads
  \[ Z (X,t) =
    \prod_{i \in \ZZ} \det \bigl(1 - tF \bigm| H^i_c (\overline{X}, \QQ_\ell)\bigr)^{(-1)^{i+1}}, \]
  where $\overline{X} \dfn X \times_{\Spec \FF_q} \overline{\FF}_q$ and $F$ is
  the Frobenius acting on $H^i_c (\overline{X}, \QQ_\ell)$. Substituting
  $t = q^{-n}$,
  \[ \zeta (X,n) =
    \prod_{i \in \ZZ} \det \bigl(1 - q^{-n} F \bigm| H^i_c (\overline{X}, \QQ_\ell)\bigr)^{(-1)^{i+1}}. \]
  Then, by the proof of \cite[Theorem~(3.1)]{Bayer-Neukirch-1978}, for each
  $\ell \ne p$, we obtain
  \begin{equation}
    \label{eqn:bayer-neukirch}
    |\zeta (X,n)|_\ell =
    \prod_{i \in \ZZ} |H^i_c (X, \ZZ_\ell (n))|^{(-1)^{i+1}}.
  \end{equation}

  On the other hand, for $n < 0$ we have
  \begin{equation}
    \label{eqn:zeta-X/Fq-p-part}
    |\zeta (X,n)|_p = 1.
  \end{equation}
  This fact can be justified, without assuming that $X$ is smooth or projective,
  e.g., using Kedlaya's trace formula for rigid cohomology
  \cite[p.\,1446]{Kedlaya-2006}, which gives
  \[ Z (X,t) = \prod_i P_i (t)^{(-1)^{i+1}},
    \quad\text{where }
    P_i (t) \in \ZZ[t] \text{ and } P_i (0) = 1. \]
  In particular, $P_i (q^{-n}) \equiv 1 \pmod{p}$.

  The product formula recovers from \eqref{eqn:bayer-neukirch} and
  \eqref{eqn:zeta-X/Fq-p-part} our special value formula
  \eqref{eqn:zeta-X/Fq-product-formula}.
\end{proof}

\begin{remark}
  The fact that $|\zeta (X,n)|_p = 1$, as observed in the argument above,
  explains why our Weil-\'{e}tale cohomology ignores the $p$-primary part in
  some sense.
\end{remark}

Let us consider a few examples to see how the special value conjecture works
over finite fields.

% Note that for a general arithmetic scheme $X$, calculating the motivic
% cohomology $H^i (X_\et, \ZZ^c (n))$ (and hence our Weil-\'{e}tale cohomology
% $H^i_\Wc (X, \ZZ(n))$) is by no means a trivial task; even the finite generation
% is known only for certain cases (see \cite[\S 8]{Beshenov-Weil-etale-1}).
% Similarly, an explicit computation of the regulator map $Reg_{X,n}$ is highly
% nontrivial. Therefore, we give some easy examples over finite fields.

\begin{example}
  \label{example:C(X,n)-for-Spec-Fq}
  If $X = \Spec \FF_q$, then $\zeta (X,s) = \frac{1}{1 - q^{-s}}$. In this case
  for $n < 0$ we obtain
  \begin{equation}
    \label{eqn:motivic-cohomology-of-Fq}
    H^i (\Spec \FF_{q,\et}, \ZZ^c (n)) \cong
    \begin{cases}
      \ZZ/(q^{-n} - 1), & i = 1, \\
      0, & i \ne 1
    \end{cases}
  \end{equation}
  (see, for example, \cite[Example~4.2]{Geisser-2017}). Therefore, formula
  \eqref{eqn:special-value-for-X/Fq} indeed recovers $\zeta (X,n)$ up to sign.

  Similarly, if we replace $\Spec \FF_q$ with $\Spec \FF_{q^m}$, considered as a
  variety over $\FF_q$, then
  $\zeta (\Spec \FF_{q^m},s) = \zeta (\Spec \FF_q, ms)$, and
  \eqref{eqn:motivic-cohomology-of-Fq} also changes accordingly.
\end{example}

\begin{example}
  Consider $X = \PP^1_{\FF_q}/(0\sim 1)$, or equivalently, a nodal cubic.
  The zeta function is $\zeta (X,s) = \frac{1}{1 - q^{1-s}}$.  We can calculate
  the groups $H^i (X_\et, \ZZ^c (n))$ using the blowup square
  \[ \begin{tikzcd}
      \Spec \FF_q \sqcup \Spec \FF_q \ar{r}\ar{d}\tikzpb & \PP^1_{\FF_q} \ar{d} \\
      \Spec \FF_q \ar{r} & X
    \end{tikzcd} \]
  This is similar to \cite[\S 8, Example~2]{Geisser-2006}. Geisser uses the
  eh-topology and long exact sequences associated to abstract blowup squares
  \cite[Proposition~3.2]{Geisser-2006}. In our case the same reasoning works,
  because by \cite[Theorem~I]{Beshenov-Weil-etale-1}, one has
  $H^i (X_\et, \ZZ^c (n)) \cong \Hom (H^{2-i}_c (X_\et, \ZZ (n)),\QQ/\ZZ)$,
  where $\ZZ (n) = \varinjlim_{p\nmid m} \mu_m^{\otimes n} [-1]$, and \'{e}tale
  cohomology and eh-cohomology coincide for such sheaves by
  \cite[Theorem~3.6]{Geisser-2006}.

  Using the projective bundle formula, we calculate from
  \eqref{eqn:motivic-cohomology-of-Fq}
  \[ H^i (\PP^1_{\FF_q,\et}, \ZZ^c (n)) \cong \begin{cases}
      \ZZ/(q^{1-n} - 1), & i = -1, \\
      \ZZ/(q^{-n} - 1), & i = +1, \\
      0, & i \ne \pm 1.
    \end{cases} \]
  By the argument from \cite[\S 8, Example~2]{Geisser-2006}, the short exact
  sequences
  \[ 0 \to H^i (\PP^1_{\FF_q,\et}, \ZZ^c (n)) \to
    H^i (X_\et, \ZZ^c (n)) \to
    H^{i+1} ((\Spec \FF_q)_\et, \ZZ^c (n)) \to 0 \]
  give
  \[ H^i (X_\et, \ZZ^c (n)) \cong \begin{cases}
      \ZZ/(q^{1-n} - 1), & i = -1, \\
      \ZZ/(q^{-n} - 1), & i = 0,1, \\
      0, & \text{otherwise}.
    \end{cases} \]
  The formula \eqref{eqn:special-value-for-X/Fq} gives the correct value
  $\zeta (X,n)$.
\end{example}

\begin{example}
  In general, if $X/\FF_q$ is a curve, then Conjecture~$\mathbf{L}^c (X_\et,n)$
  holds; see for example \cite[Proposition~4.3]{Geisser-2017}. The cohomology
  $H^i (X_\et, \ZZ^c(n))$ is concentrated in degrees $-1, 0, +1$ by duality
  \cite[Theorem~I]{Beshenov-Weil-etale-1} and the reasons of cohomological
  dimension, and the special value formula is
  \[ \zeta^* (X,n) =
    \pm \frac{|H^0 (X_\et, \ZZ^c (n))|}{|H^{-1} (X_\et, \ZZ^c (n))|\cdot |H^1 (X_\et, \ZZ^c (n))|}. \]
\end{example}

%%%%%%%%%%%%%%%%%%%%%%%%%%%%%%%%%%%%%%%%%%%%%%%%%%%%%%%%%%%%%%%%%%%%%%%%%%%%%%%%

\section{Compatibility with operations on schemes}
\label{sec:compatibility-with-operations}

The following basic properties follow from the definition of $\zeta (X,s)$
(formula \eqref{eqn:Euler-product-for-zeta}).

\begin{enumerate}
\item[1)] \textbf{Disjoint unions}: if $X = \coprod_{1 \le i \le r} X_i$ is a
  finite disjoint union of arithmetic schemes, then
  \begin{equation}
    \label{eqn:zeta-function-for-disjoint-unions}
    \zeta (X,s) = \prod_{1 \le i \le r} \zeta (X_i,s).
  \end{equation}
  In particular,
  \begin{align*}
    \ord_{s=n} \zeta (X,s) & = \sum_{1 \le i \le r} \ord_{s=n} \zeta (X_i,s), \\
    \zeta^* (X,n) & = \prod_{1 \le i \le r} \zeta^* (X_i,n).
  \end{align*}

\item[2)] \textbf{Closed-open decompositions}: if $Z \subset X$ is a closed
  subscheme and $U = X\setminus Z$ is its open complement, then we say that we
  have a \textbf{closed-open decomposition} and write
  $Z \not\hookrightarrow X \hookleftarrow U$. In this case
  \begin{equation}
    \label{eqn:zeta-function-for-closed-open-decompositions}
    \zeta (X,s) = \zeta (Z,s) \cdot \zeta (U,s).
  \end{equation}
  In particular,
  \begin{align*}
    \ord_{s=n} \zeta (X,s) & = \ord_{s=n} \zeta (Z,s) + \ord_{s=n} \zeta (U,s), \\
    \zeta^* (X,n) & = \zeta^* (Z,n) \cdot \zeta^* (U,n).
  \end{align*}

\item[3)] \textbf{Affine bundles}: for any $r \ge 0$ the zeta function of the
  relative affine space $\AA^r_X = \AA^r_\ZZ \times X$ satisfies
  \begin{equation}
    \label{eqn:zeta-function-for-affine-space}
    \zeta (\AA^r_X, s) = \zeta (X, s-r).
  \end{equation}
  In particular,
  \begin{align*}
    \ord_{s=n} \zeta (\AA^r_X, s) & = \ord_{s=n-r} \zeta (X, s), \\
    \zeta^* (\AA^r_X, n) & = \zeta^* (X, n-r).
  \end{align*}
\end{enumerate}

This suggests that Conjectures~$\mathbf{VO} (X,n)$ and $\mathbf{C} (X,n)$ should
also satisfy the corresponding compatibilities. We verify in this section that
this is indeed the case.

\begin{lemma}
  \label{lemma:compatibility-of-Lc(X,n)}
  Let $n < 0$.

  \begin{enumerate}
  \item[1)] If $X = \coprod_{1 \le i \le r} X_i$ is a finite disjoint union of
    arithmetic schemes, then
    $$\mathbf{L}^c (X_\et,n) \iff \mathbf{L}^c (X_{i,\et},n)\text{ for all }i.$$

  \item[2)] For a closed-open decomposition
    $Z \not\hookrightarrow X \hookleftarrow U$, if two of the three conjectures
    \[ \mathbf{L}^c (X_\et,n), \quad
      \mathbf{L}^c (Z_\et,n), \quad
      \mathbf{L}^c (U_\et, n) \]
    are true, then the third is also true.

  \item[3)] For an arithmetic scheme $X$ and any $r \ge 0$, one has
    $$\mathbf{L}^c (\AA^r_{X,\et}, n) \iff \mathbf{L}^c (X_\et, n-r).$$
  \end{enumerate}

  \begin{proof}
    See the proof of \cite[Proposition~5.10]{Morin-2014}.
  \end{proof}
\end{lemma}

\begin{lemma}
  \label{lemma:compatibility-of-B(X,n)}
  Let $n < 0$.

  \begin{enumerate}
  \item[1)] If $X = \coprod_{1 \le i \le r} X_i$ is a finite disjoint union of
    arithmetic schemes, then
    \begin{multline*}
      Reg_{X,n} = \bigoplus_{1 \le i \le r} Reg_{X_i,n}\colon\\
      \bigoplus_{1 \le i \le r} R\Gamma (X_{i,\et}, \RR^c (n)) \to
      \bigoplus_{i \le i \le r} R\Gamma_\BM (G_\RR, X_i (\CC), \RR (n)) [1].
    \end{multline*}
    In particular,
    $$\mathbf{B} (X,n) \iff \mathbf{B} (X_i,n)\text{ for all }i.$$

  \item[2)] For a closed-open decomposition of arithmetic schemes
    $Z \not\hookrightarrow X \hookleftarrow U$, the corresponding regulators
    give a morphism of distinguished triangles

    \[ \begin{tikzcd}[column sep=4em]
        R\Gamma (Z_\et, \RR^c (n)) \ar{d}\ar{r}{Reg_{Z,n}} & R\Gamma_\BM (G_\RR, Z (\CC), \RR (n)) [1] \ar{d} \\
        R\Gamma (X_\et, \RR^c (n)) \ar{d}\ar{r}{Reg_{X,n}} & R\Gamma_\BM (G_\RR, X (\CC), \RR (n)) [1] \ar{d} \\
        R\Gamma (U_\et, \RR^c (n)) \ar{d}\ar{r}{Reg_{U,n}} & R\Gamma_\BM (G_\RR, U (\CC), \RR (n)) [1] \ar{d} \\
        R\Gamma (Z_\et, \RR^c (n)) [1]\ar{r}{Reg_{Z,n} [1]} & R\Gamma_\BM (G_\RR, Z (\CC), \RR (n)) [2]
      \end{tikzcd} \]
    In particular, if two of the three conjectures
    \[ \mathbf{B} (X,n), \quad
      \mathbf{B} (Z,n), \quad
      \mathbf{B} (U,n) \]
    are true, then the third is also true.

  \item[3)] For any $r \ge 0$, the diagram
    \[ \begin{tikzcd}
        R\Gamma (X_\et, \RR^c (n-r)) [2r] \ar{d}{Reg_{X,n-r}}\ar{r}{\cong} & R\Gamma (\AA^r_{X,\et}, \RR^c (n))\ar{d}{Reg_{\AA^r_X,n}} \\
        R\Gamma_\BM (G_\RR, X (\CC), \RR (n-r)) [2r] \ar{r}{\cong} & R\Gamma_\BM (G_\RR, \AA^r_X (\CC), \RR (n))
      \end{tikzcd} \]
    commutes. In particular, one has
    $$\mathbf{B} (\AA^r_X, n) \iff \mathbf{B} (X, n-r).$$
  \end{enumerate}

  \begin{proof}
    Part 1) is clear because all cohomologies that enter the definition of
    $Reg_{X,n}$ decompose into direct sums over $i = 1,\ldots r$.  Parts 2) and
    3) boil down to the corresponding functoriality properties for the KLM
    morphism \eqref{eqn:KLM-morphism-1}, namely that it commutes with proper
    pushforwards and flat pullbacks by
    \cite[Lemma~3 and~4]{Weisschuh-2017}. For closed-open decompositions, the
    distinguished triangle
    \[ R\Gamma (Z_\et, \RR^c (n)) \to R\Gamma (X_\et, \RR^c (n)) \to
      R\Gamma (U_\et, \RR^c (n)) \to R\Gamma (Z_\et, \RR^c (n)) [1] \]
    comes exactly from the proper pushforward along $Z \hookrightarrow X$ and
    flat pullback along $U \hookrightarrow X$ (see
    \cite[Corollary~7.2]{Geisser-2010} and \cite[\S 3]{Bloch-1986}). Similarly,
    the quasi-isomorphism
    $R\Gamma (X_\et, \RR^c (n-r)) [2r] \cong R\Gamma (\AA^r_{X,\et}, \RR^c (n))$
    results from the flat pullback along $p\colon \AA^r_X \to X$.
  \end{proof}
\end{lemma}

\begin{proposition}
  \label{prop:compatibility-of-VO(X,n)}
  For each arithmetic scheme $X$ below and $n < 0$, assume
  $\mathbf{L}^c (X_\et,n)$, $\mathbf{B} (X,n)$, and the meromorphic continuation
  of $\zeta (X,s)$ around $s = n$.

  \begin{enumerate}
  \item[1)] If $X = \coprod_{1 \le i \le r} X_i$ is a finite disjoint union of
    arithmetic schemes, then
    $$\mathbf{VO} (X,n) \iff \mathbf{VO} (X_i,n)\text{ for all }i.$$

  \item[2)] For a closed-open decomposition
    $Z \not\hookrightarrow X \hookleftarrow U$,
    if two of the three conjectures
    \[ \mathbf{VO} (X,n), \quad
      \mathbf{VO} (Z,n), \quad
      \mathbf{VO} (U,n) \]
    are true, then the third is also true.

  \item[3)] For any $r \ge 0$, one has
    $$\mathbf{VO} (\AA^r_X, n) \iff \mathbf{VO} (X, n-r).$$
  \end{enumerate}

  \begin{proof}
    We have already observed in Proposition~\ref{prop:VO(X,n)-assuming-B(X,n)}
    that under Conjecture~$\mathbf{B} (X,n)$ we can rewrite $\mathbf{VO} (X,n)$
    as
    $$\ord_{s=n} \zeta (X,s) = \chi (R\Gamma_c (G_\RR, X(\CC), \RR (n))).$$

    In part 1), we have
    $$\ord_{s=n} \zeta (X,s) = \sum_{1 \le i \le r} \ord_{s=n} \zeta (X_i,s),$$
    and for the corresponding $G_\RR$-equivariant cohomology,
    \[ R\Gamma_c (G_\RR, X(\CC), \RR (n)) =
      \bigoplus_{1 \le i \le r} R\Gamma_c (G_\RR, X(\CC), \RR (n)). \]
    The statement follows from the additivity of the Euler characteristic:
    \[ \begin{tikzcd}[column sep=5em]
        \ord_{s=n} \zeta (X,s) \ar[equals]{r}{\mathbf{VO} (X,n)}\ar[equals]{d} & \chi (R\Gamma_c (G_\RR, X(\CC), \RR (n))) \ar[equals]{d} \\
        \sum\limits_{1 \le i \le r} \ord_{s=n} \zeta (X_i,s) \ar[equals]{r}{\forall i \mathbf{VO} (X_i,n)} & \sum\limits_{1 \le i \le r} \chi (R\Gamma_c (G_\RR, X_i (\CC), \RR (n)))
      \end{tikzcd} \]

    Similarly in part 2), we can consider the distinguished triangle
    \begin{multline*}
      R\Gamma_c (G_\RR, U (\CC), \RR (n)) \to
      R\Gamma_c (G_\RR, X (\CC), \RR (n)) \to
      R\Gamma_c (G_\RR, Z (\CC), \RR (n)) \\
      \to R\Gamma_c (G_\RR, U (\CC), \RR (n)) [1]
  \end{multline*}
  and the additivity of the Euler characteristic gives
  \[ \begin{tikzcd}[column sep=4em]
      \ord_{s=n} \zeta (X,s) \ar[equals]{r}{\mathbf{VO} (X,n)}\ar[equals]{d} & \chi (R\Gamma_c (G_\RR, X(\CC), \RR (n))) \ar[equals]{d} \\
      \ord_{s=n} \zeta (Z,s) \ar[equals]{r}{\mathbf{VO} (Z,n)} & \chi (R\Gamma_c (G_\RR, Z (\CC), \RR (n))) \\[-2em]
      + & + \\[-2em]
      \ord_{s=n} \zeta (U,s) \ar[equals]{r}{\mathbf{VO} (U,n)} & \chi (R\Gamma_c (G_\RR, U (\CC), \RR (n)))
    \end{tikzcd} \]

  Finally, in part 3), we assume for simplicity that $X_\CC$ is connected of
  dimension $d_\CC$. Then the Poincar\'{e} duality and homotopy invariance of
  the usual cohomology without compact support give us
  \begin{multline*}
    R\Gamma_c (G_\RR, \AA^r (\CC) \times X (\CC), \RR (n)) \\
    \stackrel{\text{P.D.}}{\cong}
    \RHom (R\Gamma (G_\RR, \AA^r (\CC) \times X (\CC), \RR (d_\CC + r - n)), \RR) [-2d_\CC - 2r] \\
    \stackrel{\text{H.I.}}{\cong}
    \RHom (R\Gamma (G_\RR, X (\CC), \RR (d_\CC + r - n)), \RR) [-2d_\CC - 2r] \\
    \stackrel{\text{P.D.}}{\cong}
    R\Gamma_c (G_\RR, X (\CC), \RR (n - r)) [-2r].
  \end{multline*}
  The twist $[-2r]$ is even and therefore has no effect on the Euler
  characteristic, so that we obtain
  \[ \begin{tikzcd}[column sep=4em]
      \ord_{s=n} \zeta (\AA^r_X,s) \ar[equals]{r}{\mathbf{VO} (\AA^r_X,n)}\ar[equals]{d} & \chi (R\Gamma_c (G_\RR, \AA^r (\CC) \times X(\CC), \RR (n))) \ar[equals]{d} \\
      \ord_{s=n-r} \zeta (X,s) \ar[equals]{r}{\mathbf{VO} (X,n-r)} & \chi (R\Gamma_c (G_\RR, X (\CC), \RR (n-r)))
    \end{tikzcd} \]
  \end{proof}
\end{proposition}

% \begin{remark}
%   Recall that the formula that appears in the original statement of
%   $\mathbf{VO} (X,n)$ reads
%   \begin{equation}
%     \label{eqn:VO(X,n)-original-formula}
%     \ord_{s=n} \zeta (X,s) = \chi' (R\Gamma_\Wc (X,\ZZ(n))) \dfn
%     \sum_{i\in \ZZ} (-1)^i\cdot i \cdot \rk_\ZZ H^i_\Wc (X,\ZZ(n)).
%   \end{equation}
%   Conjecture~$\mathbf{B} (X,n)$ in the above argument is needed to rewrite
%   this in terms of the usual Euler characteristic. We used
%   $\chi (R\Gamma_c (G_\RR, X (\CC), \RR (n)))$, but we could do the same with
%   $\chi (\RHom (R\Gamma (X_\et, \ZZ^c (n)), \RR) [1])$.

%   The least interesting part 1) of the previous proposition could be proved
%   directly from \eqref{eqn:VO(X,n)-original-formula}, since
%   $H^i_\Wc (X,\ZZ(n)) = \bigoplus_j H^i_\Wc (X_j,\ZZ(n))$. Parts 2) and 3) would
%   be problematic to prove directly from \eqref{eqn:VO(X,n)-original-formula}
%   without assuming $\mathbf{B} (X,n)$, since the secondary Euler characteristic
%   $\chi' (-)$ does not behave as the usual Euler characteristic $\chi (-)$.
%   In particular, it is not additive for distinguished triangles.
% \end{remark}

Our next goal is to prove similar compatibilities for
Conjecture~$\mathbf{C} (X,n)$, as was just done for $\mathbf{VO} (X,n)$.
We split the proof into three technical lemmas
\ref{lemma:lambda-and-disjoint-unions},
\ref{lemma:lambda-and-closed-open-decompositions},
\ref{lemma:lambda-and-affine-bundles}, each for the corresponding compatibility.

\begin{lemma}
  \label{lemma:lambda-and-disjoint-unions}
  Let $n < 0$ and let $X = \coprod_{1 \le i \le r} X_i$ be a finite disjoint
  union of arithmetic schemes. Assume $\mathbf{L}^c (X_\et,n)$ and
  $\mathbf{B} (X,n)$. Then there is a quasi-isomorphism of complexes
  \begin{equation}
    \label{eqn:RGamma-Wc-and-disjoint-unions}
    \bigoplus_{1 \le i \le r} R\Gamma_\Wc (X_i, \ZZ(n)) \cong
    R\Gamma_\Wc (X, \ZZ(n)),
  \end{equation}
  which after taking the determinants gives a commutative diagram

  \begin{equation}
    \label{eqn:lambda-and-disjoint-unions}
    \begin{tikzcd}
      \RR \otimes_\RR \cdots \otimes_\RR \RR\ar{d}{\lambda_{X_1,n}\otimes\cdots\otimes\lambda_{X_r,n}}[swap]{\cong} \ar{r}{x_1\otimes\cdots\otimes x_r \mapsto x_1\cdots x_r}[swap]{\cong} & \RR \ar{d}{\lambda_{X,n}}[swap]{\cong} \\
      \bigotimes\limits_{1 \le i \le r} (\det_\ZZ R\Gamma_\Wc (X_i, \ZZ(n))) \otimes \RR \ar{r}{\cong} & (\det_\ZZ R\Gamma_\Wc (X_i, \ZZ(n))) \otimes \RR
    \end{tikzcd}
  \end{equation}

  \begin{proof}
    For $X = \coprod_{1 \le i \le r} X_i$, all cohomologies in our construction
    of $R\Gamma_\Wc (X, \ZZ(n))$ in \cite{Beshenov-Weil-etale-1} decompose into
    the corresponding direct sum over $i = 1,\ldots,r$, and
    \eqref{eqn:RGamma-Wc-and-disjoint-unions} follows.

    After tensoring with $\RR$, we obtain a commutative diagram

    \[ \begin{tikzcd}[column sep=1.25em]
        \bigoplus_i \left(\!\!\!\begin{array}{c} R\Gamma_c (G_\RR, X_i (\CC), \RR (n)) [-2] \\ \oplus \\ R\Gamma_c (G_\RR, X_i (\CC), \RR (n)) [-1] \end{array}\!\!\!\right) \ar{d}{\bigoplus_i Reg_{X_i,n}^\vee [-1] \oplus id}[swap]{\cong} \ar{r}{\cong} & \begin{array}{c} R\Gamma_c (G_\RR, X (\CC), \RR (n)) [-2] \\ \oplus \\ R\Gamma_c (G_\RR, X (\CC), \RR (n)) [-1] \end{array} \ar{d}{Reg_{X,n}^\vee [-1] \oplus id}[swap]{\cong} \\
        \bigoplus_i \left(\!\!\!\begin{array}{c} \RHom (R\Gamma (X_{i,\et}, \ZZ^c (n)), \RR) [-1] \\ \oplus \\ R\Gamma_c (G_\RR, X_i (\CC), \RR (n)) [-1] \end{array}\!\!\!\right) \ar{d}{\text{split}}[swap]{\cong} \ar{r}{\cong} & \begin{array}{c} \RHom (R\Gamma (X_\et, \ZZ^c (n)), \RR) [-1] \\ \oplus \\ R\Gamma_c (G_\RR, X (\CC), \RR (n)) [-1] \end{array} \ar{d}{\text{split}}[swap]{\cong} \\
        \bigoplus_i R\Gamma_\Wc (X_i, \RR (n)) \ar{r}{\cong} & R\Gamma_\Wc (X, \RR (n))
      \end{tikzcd} \]
    Taking the determinants, we obtain \eqref{eqn:lambda-and-disjoint-unions}.
  \end{proof}
\end{lemma}

\begin{lemma}
  \label{lemma:lambda-and-closed-open-decompositions}
  Let $n < 0$ and let $Z \not\hookrightarrow X \hookleftarrow U$ be a
  closed-open decomposition of arithmetic schemes, such that the conjectures
  \begin{gather*}
    \mathbf{L}^c (U_\et,n), ~ \mathbf{L}^c (X_\et,n), ~ \mathbf{L}^c (Z_\et,n),\\
    \mathbf{B} (U,n), ~ \mathbf{B} (X,n), ~ \mathbf{B} (Z_\et,n)
  \end{gather*}
  hold (it suffices to assume two of the three conjectures thanks to Lemmas
  \ref{lemma:compatibility-of-Lc(X,n)} and \ref{lemma:compatibility-of-B(X,n)}).
  Then there is an isomorphism of determinants
  \begin{equation}
    \label{eqn:isomorphism-of-det-RGamma-Wc-for-closed-open-decompositions}
    \det_\ZZ R\Gamma_\Wc (U, \ZZ(n)) \otimes
    \det_\ZZ R\Gamma_\Wc (Z, \ZZ(n)) \cong
    \det_\ZZ R\Gamma_\Wc (X, \ZZ(n))
  \end{equation}
  making the following diagram commute up to signs:
  \begin{equation}
    \label{eqn:lambda-and-closed-open-decompositions}
    \begin{tikzcd}
      \RR \otimes_\RR \RR \ar{r}{x\otimes y \mapsto xy}\ar{d}{\lambda_{U,n} \otimes \lambda_{Z,n}}[swap]{\cong} & \RR\ar{d}{\lambda_{X,n}}[swap]{\cong} \\
      \begin{array}{c} (\det_\ZZ R\Gamma_\Wc (U, \ZZ(n)))\otimes \RR \\ \otimes_\RR \\ (\det_\ZZ R\Gamma_\Wc (Z, \ZZ(n)))\otimes \RR \end{array} \ar{r}{\cong} & (\det_\ZZ R\Gamma_\Wc (X, \ZZ(n))) \otimes \RR
    \end{tikzcd}
  \end{equation}

  \begin{proof}
    A closed-open decomposition $Z \not\hookrightarrow X \hookleftarrow U$
    induces the distinguished triangles
    \[ \begin{tikzcd}[row sep=0pt,column sep=1em,font=\small]
        R\Gamma (Z_\et, \ZZ^c (n)) \ar{r} & R\Gamma (X_\et, \ZZ^c (n)) \ar{r} & R\Gamma (U_\et, \ZZ^c (n)) \ar{r} & {[1]} \\
        R\Gamma_c (U_\et, \ZZ (n)) \ar{r} & R\Gamma_c (X_\et, \ZZ (n)) \ar{r} & R\Gamma_c (Z_\et, \ZZ (n)) \ar{r} & {[1]} \\
        R\Gamma_c (G_\RR, U (\CC), \RR (n)) \ar{r} & R\Gamma_c (G_\RR, X (\CC), \RR (n)) \ar{r} & R\Gamma_c (G_\RR, Z (\CC), \RR (n)) \ar{r} & {[1]}
      \end{tikzcd} \]
    The first triangle is \cite[Corollary~7.2]{Geisser-2010} and it means that
    $R\Gamma (-, \ZZ^c (n))$ behaves like Borel--Moore homology. The following
    two are the usual triangles for cohomology with compact support. These fit
    together in a commutative diagram shown in
    Figure~\ref{fig:RGamma-Wc-and-closed-open-decompositions} below
    (p.~\pageref{fig:RGamma-Wc-and-closed-open-decompositions}). Figure~\ref{fig:RGamma-Wc-and-closed-open-decompositions-otimes-Q}
    on p.~\pageref{fig:RGamma-Wc-and-closed-open-decompositions-otimes-Q} shows
    the same diagram tensored with $\RR$.

    In this diagram we start from the morphism of triangles
    $(\alpha_{U,n}, \alpha_{X,n}, \alpha_{Z,n})$ and then take the corresponding
    cones $R\Gamma_\fg (-, \ZZ(n))$. By
    \cite[Proposition~5.6]{Beshenov-Weil-etale-1}, these cones are defined up to
    a \emph{unique} isomorphism in the derived category $\mathbf{D} (\ZZ)$, and
    the same argument shows that the induced morphisms of complexes
    \begin{equation}
      \label{eqn:triangle-RGamma-fg}
      R\Gamma_\fg (U, \ZZ(n)) \to
      R\Gamma_\fg (X, \ZZ(n)) \to
      R\Gamma_\fg (Z, \ZZ(n)) \to
      R\Gamma_\fg (U, \ZZ(n)) [1]
    \end{equation}
    are also well-defined
    (see \cite[Corollary~A.3]{Beshenov-Weil-etale-1}). A priori,
    \eqref{eqn:triangle-RGamma-fg} need not be a distinguished triangle, but we
    claim that it induces a long exact sequence in cohomology.

    To this end, note that tensoring the diagram with $\ZZ/m\ZZ$ gives us an
    isomorphism
    \[ \begin{tikzcd}[column sep=1em,font=\small]
        R\Gamma_c (U_\et, \ZZ/m\ZZ (n)) \ar{r}\ar{d}{\cong} & R\Gamma_c (X_\et, \ZZ/m\ZZ (n)) \ar{r}\ar{d}{\cong} & R\Gamma_c (Z_\et, \ZZ/m\ZZ (n)) \ar{r}\ar{d}{\cong} & {[1]}\ar{d}{\cong} \\
        \begin{array}{c} R\Gamma_\fg (U, \ZZ (n)) \\ \otimes^\mathbf{L} \\ \ZZ/m\ZZ \end{array} \ar{r} & \begin{array}{c} R\Gamma_\fg (X, \ZZ (n)) \\ \otimes^\mathbf{L} \\ \ZZ/m\ZZ \end{array} \ar{r} & \begin{array}{c} R\Gamma_\fg (Z, \ZZ (n)) \\ \otimes^\mathbf{L} \\ \ZZ/m\ZZ \end{array} \ar{r} & {[1]}
      \end{tikzcd} \]
    More generally, for each prime $p$ we can take the corresponding derived
    $p$-adic completions (see \cite{Bhatt-Scholze-2015} and
    \cite[Tag~091N]{Stacks-project})
    \[ R\Gamma_\fg (-, \ZZ(n))^\wedge_p \dfn
      R\varprojlim_k (R\Gamma_\fg (-, \ZZ(n)) \otimes^\mathbf{L} \ZZ/p^k\ZZ), \]
    which give us a distinguished triangle for each prime $p$
    \[ R\Gamma_\fg (U, \ZZ(n))^\wedge_p \to
      R\Gamma_\fg (X, \ZZ(n))^\wedge_p \to
      R\Gamma_\fg (Z, \ZZ(n))^\wedge_p \to
      R\Gamma_\fg (U, \ZZ(n))^\wedge_p [1]. \]
    At the level of cohomology, there are natural isomorphisms
    \cite[Tag~0A06]{Stacks-project}
    \[ H^i (R\Gamma_\fg (-, \ZZ(n))^\wedge_p) \cong
      H^i_\fg (-, \ZZ(n)) \otimes \ZZ_p. \]
    In particular, for each $p$ there is a long exact sequence of cohomology
    groups
    \begin{multline*}
      \cdots \to H^i_\fg (U, \ZZ(n)) \otimes \ZZ_p \to
      H^i_\fg (X, \ZZ(n)) \otimes \ZZ_p \to
      H^i_\fg (Z, \ZZ(n)) \otimes \ZZ_p \\
      \to H^{i+1}_\fg (U, \ZZ(n)) \otimes \ZZ_p \to \cdots
    \end{multline*}
    induced by \eqref{eqn:triangle-RGamma-fg}. By finite generation of
    $H^i_\fg (-, \ZZ(n))$ and flatness of $\ZZ_p$ this implies that the sequence
    \begin{equation}
      \label{eqn:RGamma-fg-long-exact-sequence}
      \cdots \to H^i_\fg (U, \ZZ(n)) \to
      H^i_\fg (X, \ZZ(n)) \to
      H^i_\fg (Z, \ZZ(n)) \to
      H^{i+1}_\fg (U, \ZZ(n)) \to \cdots
    \end{equation}
    is exact.

    Now we consider the diagram
    \[ \begin{tikzcd}[column sep=0.75em,font=\small]
        \tau_{\le m} R\Gamma_c (G_\RR, U (\CC), \ZZ (n))[-1] \ar{r}\ar{d} & R\Gamma_\Wc (U, \ZZ (n))\ar{d}\ar{r} & \tau_{\le m} R\Gamma_\fg (U,\ZZ(n)) \ar{d}\ar{r} & {[0]} \ar{d} \\
        \tau_{\le m} R\Gamma_c (G_\RR, X (\CC), \ZZ (n))[-1] \ar{r}\ar{d} & R\Gamma_\Wc (X, \ZZ (n))\ar{d}\ar{r} & \tau_{\le m} R\Gamma_\fg (X,\ZZ(n)) \ar{d}\ar{r} & {[0]} \ar{d} \\
        \tau_{\le m} R\Gamma_c (G_\RR, Z (\CC), \ZZ (n))[-1] \ar{r}\ar{d} & R\Gamma_\Wc (Z, \ZZ (n))\ar{d}\ar{r} & \tau_{\le m} R\Gamma_\fg (Z,\ZZ(n)) \ar{d}\ar{r} & {[0]} \ar{d} \\
        \tau_{\le m} R\Gamma_c (G_\RR, U (\CC), \ZZ (n)) \ar{r} & R\Gamma_\Wc (U, \ZZ (n)) [1] \ar{r} & \tau_{\le m} R\Gamma_\fg (U,\ZZ(n)) [1] \ar{r} & {[1]}
      \end{tikzcd} \]
    Here we took truncations for $m$ big enough, as in the proof of
    Lemma~\ref{lemma:determinant-of-RGamma-Wc-well-defined}. There are canonical
    isomorphisms
    \begin{align*}
      \notag \det_\ZZ R\Gamma_\Wc (U, \ZZ(n)) & \cong \begin{array}{c} \det_\ZZ (\tau_{\le m} R\Gamma_c (G_\RR, U (\CC), \ZZ (n)) [-1]) \\ \otimes \\ \det_\ZZ (\tau_{\le m} R\Gamma_\fg (U, \ZZ(n))), \end{array} \\
      \\
      \notag \det_\ZZ R\Gamma_\Wc (X, \ZZ(n)) & \cong \begin{array}{c} \det_\ZZ (\tau_{\le m} R\Gamma_c (G_\RR, X (\CC), \ZZ (n)) [-1]) \\ \otimes \\ \det_\ZZ (\tau_{\le m} R\Gamma_\fg (X, \ZZ(n))), \end{array} \\
      \\
      \notag \det_\ZZ R\Gamma_\Wc (Z, \ZZ(n)) & \cong \begin{array}{c} \det_\ZZ (\tau_{\le m} R\Gamma_c (G_\RR, Z (\CC), \ZZ (n)) [-1]) \\ \otimes \\ \det_\ZZ (\tau_{\le m} R\Gamma_\fg (Z, \ZZ(n))), \end{array} \\
      \\
      \det_\ZZ (\tau_{\le m} R\Gamma_c (G_\RR, X (\CC), \ZZ(n))) & \cong \begin{array}{c} \det_\ZZ (\tau_{\le m} R\Gamma_c (G_\RR, U (\CC), \ZZ (n))) \\ \otimes \\ \det_\ZZ (\tau_{\le m} R\Gamma_c (G_\RR, Z (\CC), \ZZ(n))), \end{array} \\
      \\
      \det_\ZZ (\tau_{\le m} R\Gamma_\fg (X, \ZZ(n))) & \cong \begin{array}{c} \det_\ZZ (\tau_{\le m} R\Gamma_\fg (U, \ZZ (n))) \\ \otimes \\ \det_\ZZ (\tau_{\le m} R\Gamma_\fg (Z, \ZZ(n))).\end{array}
    \end{align*}
    The first four isomorphisms arise from the corresponding distinguished
    triangles, while the last isomorphism comes from the long exact sequence
    \eqref{eqn:RGamma-fg-long-exact-sequence}, which gives an isomorphism
    \begin{equation}\small
      \bigotimes_{i \le m}
      \Bigl(\det_\ZZ H^i_\fg (U, \ZZ(n))^{(-1)^i} \otimes
      \det_\ZZ H^i_\fg (X, \ZZ(n))^{(-1)^{i+1}} \otimes
      \det_\ZZ H^i_\fg (Z, \ZZ(n))^{(-1)^i}\Bigr) \cong \ZZ.
    \end{equation}
    We can rearrange the terms
    (at the expense of introducing a $\pm 1$ sign),
    to obtain
    \begin{multline*}
      \det_\ZZ (\tau_{\le m} R\Gamma_\fg (X, \ZZ(n))) \cong
      \bigotimes_{i \le m} \det_\ZZ H^i_\fg (X, \ZZ(n)) \cong \\
      \bigotimes_{i \le m} \det_\ZZ H^i_\fg (U, \ZZ(n)) \otimes
      \bigotimes_{i \le m} \det_\ZZ H^i_\fg (Z, \ZZ(n)) \cong \\
      \det_\ZZ (\tau_{\le m} R\Gamma_\fg (U, \ZZ(n))) \otimes
      \det_\ZZ (\tau_{\le m} R\Gamma_\fg (Z, \ZZ(n))).
    \end{multline*}

    All this gives us the desired isomorphism of integral determinants
    \eqref{eqn:isomorphism-of-det-RGamma-Wc-for-closed-open-decompositions}.

    Let us now consider the diagram with distinguished rows in
    Figure~\ref{fig:Regulators-and-closed-open-decompositions}
    (p.~\pageref{fig:Regulators-and-closed-open-decompositions}).
    Here the three squares with the regulators involved commute thanks to
    Lemma~\ref{lemma:compatibility-of-B(X,n)}. Taking the determinants, we
    obtain \eqref{eqn:lambda-and-closed-open-decompositions}, by the
    compatibility with distinguished triangles.
  \end{proof}
\end{lemma}

\begin{remark}
  Morally, we expect that a closed-open decomposition induces a distinguished
  triangle of the form
  \begin{equation}
    \label{eqn:closed-open-decompositions-hypothetical-RGamma-triangle}
    R\Gamma_\Wc (U, \ZZ(n)) \to
    R\Gamma_\Wc (X, \ZZ(n)) \to
    R\Gamma_\Wc (Z, \ZZ(n)) \to [1].
  \end{equation}
  However, $R\Gamma_\Wc (X, \ZZ(n))$ is defined in \cite{Beshenov-Weil-etale-1}
  as a mapping fber of a morphism in $\mathbf{D} (\ZZ)$, so it is not quite
  functorial.

  We recall that in the usual derived (1-)category $\mathbf{D} (\mathcal{A})$,
  taking naively a ``cone of a morphism of distinguished triangles''
  $$\begin{tikzpicture}[ampersand replacement=\&]
    \matrix(m)[matrix of math nodes, row sep=1.5em, column sep=1.5em,
    text height=1.5ex, text depth=0.25ex]{
      A^\bullet \& B^\bullet \& C^\bullet \& A^\bullet[1] \\
      A^{\bullet\prime} \& B^{\bullet\prime} \& C^{\bullet\prime} \& A^{\bullet\prime} [1] \\
      A^{\bullet\prime\prime} \& B^{\bullet\prime\prime} \& C^{\bullet\prime\prime} \& A^{\bullet\prime\prime} [1] \\
      A^\bullet [1]  \& B^\bullet [1]  \& C^\bullet [1]  \& A^\bullet [2]\\};

    \path[->] (m-1-1) edge (m-1-2);
    \path[->] (m-1-2) edge (m-1-3);
    \path[->] (m-1-3) edge (m-1-4);

    \path[->] (m-2-1) edge (m-2-2);
    \path[->] (m-2-2) edge (m-2-3);
    \path[->] (m-2-3) edge (m-2-4);

    \path[dashed,->] (m-3-1) edge (m-3-2);
    \path[dashed,->] (m-3-2) edge (m-3-3);
    \path[dashed,->] (m-3-3) edge (m-3-4);

    \path[->] (m-4-1) edge (m-4-2);
    \path[->] (m-4-2) edge (m-4-3);
    \path[->] (m-4-3) edge (m-4-4);

    \path[->] (m-1-1) edge (m-2-1);
    \path[->] (m-1-2) edge (m-2-2);
    \path[->] (m-1-3) edge (m-2-3);
    \path[->] (m-1-4) edge (m-2-4);

    \path[dashed,->] (m-2-1) edge (m-3-1);
    \path[dashed,->] (m-2-2) edge (m-3-2);
    \path[dashed,->] (m-2-3) edge (m-3-3);
    \path[dashed,->] (m-2-4) edge (m-3-4);

    \path[dashed,->] (m-3-1) edge (m-4-1);
    \path[dashed,->] (m-3-2) edge (m-4-2);
    \path[dashed,->] (m-3-3) edge (m-4-3);
    \path[dashed,->] (m-3-4) edge (m-4-4);

    \node[font=\small] at ($(m-3-3)!.5!(m-4-4)$) {($-$)};
  \end{tikzpicture}$$
  usually \emph{does not} yield a distinguished triangle
  $A^{\bullet\prime\prime} \to B^{\bullet\prime\prime} \to C^{\bullet\prime\prime} \to A^{\bullet\prime\prime} [1]$.
  For a thorough discussion of this problem, see \cite{Neeman-1991}.

  For lack of a better definition for
  $R\Gamma_\Wc (X, \ZZ(n))$, we constructed the isomorphism
  \eqref{eqn:isomorphism-of-det-RGamma-Wc-for-closed-open-decompositions}
  ad hoc, without the hypothetical triangle
  \eqref{eqn:closed-open-decompositions-hypothetical-RGamma-triangle}.
\end{remark}

\begin{landscape}
  \begin{figure}
    \[ \begin{tikzcd}[font=\small]
        &[2em] &[-2.5em] &[-2.5em] &[-2.5em] R\Gamma_\Wc (U, \ZZ (n))\ar{dl} &[-2.5em] \\
        \RHom (R\Gamma (U_\et, \ZZ^c (n)), \QQ[-2]) \ar{r}{\alpha_{U,n}}\ar{dd} & R\Gamma_c (U_\et, \ZZ (n)) \ar{dr}{u_\infty} \ar{dd} \ar{rr} & & R\Gamma_\fg (U, \ZZ (n))\ar{dl}[swap]{i_\infty} \ar[dashed]{dd}{\exists!} \ar[crossing over]{rr} & & \cdots [-1] \ar{dd} \\
        & & R\Gamma_c (G_\RR, U (\CC), \ZZ (n)) & & R\Gamma_\Wc (X, \ZZ (n))\ar{dl} \\
        \RHom(R\Gamma (X_\et, \ZZ^c (n)), \QQ[-2]) \ar{r}{\alpha_{X,n}}\ar{dd} & R\Gamma_c (X_\et, \ZZ (n)) \ar{dr}{u_\infty} \ar{dd} \ar{rr} & & R\Gamma_\fg (X, \ZZ (n)) \ar{dl}[swap]{i_\infty} \ar[dashed]{dd}{\exists!} \ar[crossing over]{rr} & & \cdots [-1] \ar{dd} \\
        & & R\Gamma_c (G_\RR, X (\CC), \ZZ (n)) \ar[<-,near end,crossing over]{uu} & & R\Gamma_\Wc (Z, \ZZ (n))\ar{dl} \\
        \RHom(R\Gamma (Z_\et, \ZZ^c (n)), \QQ[-2]) \ar{r}{\alpha_{Z,n}}\ar{dd} & R\Gamma_c (Z_\et, \ZZ (n)) \ar{dr}{u_\infty} \ar{dd} \ar{rr} & & R\Gamma_\fg (Z, \ZZ (n)) \ar{dl}[swap]{i_\infty} \ar[dashed]{dd}{\exists!} \ar[crossing over]{rr} & & \cdots [-1] \ar{dd} \\
        & & R\Gamma_c (G_\RR, Z (\CC), \ZZ (n)) \ar[<-,near end,crossing over]{uu} & & R\Gamma_\Wc (U, \ZZ (n)) [1]\ar{dl} \\
        \RHom(R\Gamma (U_\et, \ZZ^c (n)), \QQ[-1]) \ar{r}{\alpha_{U,n} [1]} & R\Gamma_c (U_\et, \ZZ (n)) [1] \ar{dr}{u_\infty} \ar{rr} & & R\Gamma_\fg (U, \ZZ (n)) [1] \ar{dl}[swap]{i_\infty [1]} \ar{rr} & & \cdots [0] \\
        & & R\Gamma_c (G_\RR, U (\CC), \ZZ (n)) [1] \ar[<-,near end,crossing over]{uu} \\
      \end{tikzcd} \]

    \caption{Diagram induced by a closed-open decomposition
      $Z \not\hookrightarrow X \hookleftarrow U$.}
    \label{fig:RGamma-Wc-and-closed-open-decompositions}
  \end{figure}
\end{landscape}

\begin{landscape}
  \begin{figure}
    \[ \begin{tikzcd}[font=\small]
        &[2em] &[-2.5em] &[-2.5em] &[-2.5em] R\Gamma_\Wc (U, \RR (n))\ar{dl}\ar{dd} &[-2.5em] \\
        \RHom (R\Gamma (U_\et, \ZZ^c (n)), \RR[-2]) \ar{r}\ar{dd} & 0 \ar{dr} \ar{dd} \ar{rr} & & R\Gamma_\fg (U, \RR (n))\ar{dl}[swap]{0} \ar[dashed]{dd}{\exists!} \ar[crossing over,near start]{rr}{\cong} & & \RHom (R\Gamma (U_\et, \ZZ^c (n)), \RR[-1]) \ar{dd} \\
        & & R\Gamma_c (G_\RR, U (\CC), \RR (n)) & & R\Gamma_\Wc (X, \RR (n))\ar{dl}\ar{dd} \\
        \RHom (R\Gamma (X_\et, \ZZ^c (n)), \RR[-2]) \ar{r}\ar{dd} & 0 \ar{dr} \ar{dd} \ar{rr} & & R\Gamma_\fg (X, \RR (n)) \ar{dl}[swap]{0} \ar[dashed]{dd}{\exists!} \ar[crossing over,near start]{rr}{\cong} & & \RHom (R\Gamma (X_\et, \ZZ^c (n)), \RR[-1]) \ar{dd} \\
        & & R\Gamma_c (G_\RR, X (\CC), \RR (n)) \ar[<-,near end,crossing over]{uu} & & R\Gamma_\Wc (Z, \RR (n))\ar{dl}\ar{dd} \\
        \RHom (R\Gamma (Z_\et, \ZZ^c (n)), \RR[-2]) \ar{r}\ar{dd} & 0 \ar{dr} \ar{dd} \ar{rr} & & R\Gamma_\fg (Z, \RR (n)) \ar{dl}[swap]{0} \ar[dashed]{dd}{\exists!} \ar[crossing over,near start]{rr}{\cong} & & \RHom (R\Gamma (Z_\et, \ZZ^c (n)), \RR[-1]) \ar{dd} \\
        & & R\Gamma_c (G_\RR, Z (\CC), \RR (n)) \ar[<-,near end,crossing over]{uu} & & R\Gamma_\Wc (U, \RR (n)) [1]\ar{dl} \\
        \RHom (R\Gamma (U_\et, \ZZ^c (n)), \RR[-1]) \ar{r} & 0 \ar{dr} \ar{rr} & & R\Gamma_\fg (U, \RR (n)) [1] \ar{dl}[swap]{0} \ar{rr}{\cong} & & \RHom (R\Gamma (U_\et, \ZZ^c (n)), \RR) \\
        & & R\Gamma_c (G_\RR, U (\CC), \RR (n)) [1] \ar[<-,near end,crossing over]{uu} \\
      \end{tikzcd} \]

    \caption{Diagram induced by a closed-open decomposition
      $Z \not\hookrightarrow X \hookleftarrow U$, tensored with $\RR$.}
    \label{fig:RGamma-Wc-and-closed-open-decompositions-otimes-Q}
  \end{figure}
\end{landscape}

\begin{landscape}
  \begin{figure}
    \[ \begin{tikzcd}[column sep=1em,font=\small]
        \begin{array}{c} R\Gamma_c (G_\RR, U (\CC), \RR (n)) [-2] \\ \oplus \\ R\Gamma_c (G_\RR, U (\CC), \RR (n)) [-1] \end{array} \ar{d}{Reg_{U,n}^\vee [-1] \oplus id}[swap]{\cong} \ar{r} & \begin{array}{c} R\Gamma_c (G_\RR, X (\CC), \RR (n)) [-2] \\ \oplus \\ R\Gamma_c (G_\RR, X (\CC), \RR (n)) [-1] \end{array} \ar{d}{Reg_{X,n}^\vee [-1] \oplus id}[swap]{\cong} \ar{r} & \begin{array}{c} R\Gamma_c (G_\RR, Z (\CC), \RR (n)) [-2] \\ \oplus \\ R\Gamma_c (G_\RR, Z (\CC), \RR (n)) [-1] \end{array} \ar{d}{Reg_{Z,n}^\vee [-1] \oplus id}[swap]{\cong} \ar{r} & \begin{array}{c} R\Gamma_c (G_\RR, U (\CC), \RR (n)) [-1] \\ \oplus \\ R\Gamma_c (G_\RR, U (\CC), \RR (n)) \end{array} \ar{d}{Reg_{U,n}^\vee \oplus id}[swap]{\cong} \\
        \begin{array}{c} \RHom (R\Gamma (U_\et, \ZZ^c (n)), \RR) [-1] \\ \oplus \\ R\Gamma_c (G_\RR, U (\CC), \RR (n)) [-1] \end{array} \ar{d}{\text{split}}[swap]{\cong} \ar{r} & \begin{array}{c} \RHom (R\Gamma (X_\et, \ZZ^c (n)), \RR) [-1] \\ \oplus \\ R\Gamma_c (G_\RR, X (\CC), \RR (n)) [-1] \end{array} \ar{d}{\text{split}}[swap]{\cong} \ar{r} & \begin{array}{c} \RHom (R\Gamma (Z_\et, \ZZ^c (n)), \RR) [-1] \\ \oplus \\ R\Gamma_c (G_\RR, Z (\CC), \RR (n)) [-1] \end{array} \ar{d}{\text{split}}[swap]{\cong} \ar{r} & \begin{array}{c} \RHom (R\Gamma (U_\et, \ZZ^c (n)), \RR) \\ \oplus \\ R\Gamma_c (G_\RR, U (\CC), \RR (n)) \end{array} \ar{d}{\text{split}}[swap]{\cong} \\
        R\Gamma_\Wc (U, \RR (n)) \ar{r} & R\Gamma_\Wc (X, \RR (n)) \ar{r} & R\Gamma_\Wc (Z, \RR (n)) \ar{r} & R\Gamma_\Wc (U, \RR (n)) [1]
      \end{tikzcd} \]

    \caption{Diagram induced by a closed-open decomposition
      $Z \not\hookrightarrow X \hookleftarrow U$}
    \label{fig:Regulators-and-closed-open-decompositions}
  \end{figure}
\end{landscape}

\begin{lemma}
  \label{lemma:lambda-and-affine-bundles}
  For $n < 0$ and $r \ge 0$, let $X$ be an arithmetic scheme satisfying
  $\mathbf{L}^c (X_\et,n-r)$ and $\mathbf{B} (X,n-r)$. Then there is a natural
  quasi-isomorphism of complexes
  \begin{equation}
    \label{eqn:RGamma-Wc-and-affine-bundles}
    R\Gamma_\Wc (\AA^r_X, \ZZ (n)) \cong R\Gamma_\Wc (X, \ZZ (n-r)) [-2r],
  \end{equation}
  which after passing to the determinants makes the following diagram commute:
  \begin{equation}
    \label{eqn:lambda-and-affine-bundles}
    \begin{tikzcd}[column sep=0.5em]
    & \RR\ar{dl}{\cong}[swap]{\lambda_{\AA^r_X,n}}\ar{dr}{\lambda_{X,n-r}}[swap]{\cong} \\
      (\det_\ZZ R\Gamma_\Wc (\AA^r_X, \ZZ (n)))\otimes \RR \ar{rr}{\cong} & & (\det_\ZZ R\Gamma_\Wc (X, \ZZ (n-r)))\otimes \RR
    \end{tikzcd}
  \end{equation}

  \begin{proof}
    We refer to Figure~\ref{fig:RGamma-Wc-and-affine-bundles}
    (p.~\pageref{fig:RGamma-Wc-and-affine-bundles}), which shows how the flat
    morphism $p\colon \AA^r_X \to X$ induces the desired quasi-isomorphism
    \eqref{eqn:RGamma-Wc-and-affine-bundles}. It all boils down to the homotopy
    property of motivic cohomology, namely the fact that $p$ induces a
    quasi-isomorphism
    \[ p^*\colon R\Gamma (X_\et, \ZZ^c (n-r)) [2r] \xrightarrow{\cong}
      R\Gamma (\AA^r_{X,\et}, \ZZ^c (n)); \]
    see, e.g. \cite[Lemma~5.11]{Morin-2014}.
    After passing to real coefficients, we obtain the following diagram:
    \[ \begin{tikzcd}[column sep=1em, font=\small]
        \begin{array}{c} R\Gamma_c (G_\RR, \AA^r_X (\CC), \RR (n)) [-2] \\ \oplus \\ R\Gamma_c (G_\RR, \AA^r_X (\CC), \RR (n)) [-1] \end{array} \ar{d}{Reg_{\AA^r_X,n}^\vee [-1] \oplus id}[swap]{\cong} \ar{r}{\cong} & \begin{array}{c} R\Gamma_c (G_\RR, X (\CC), \RR (n-r)) [-2] [-2r] \\ \oplus \\ R\Gamma_c (G_\RR, X (\CC), \RR (n-r)) [-1] [-2r] \end{array} \ar{d}{Reg_{X,n-r}^\vee [-1] [-2r] \oplus id}[swap]{\cong} \\
        \begin{array}{c} \RHom (R\Gamma (\AA^r_{X,\et}, \ZZ^c (n)), \RR) [-1] \\ \oplus \\ R\Gamma_c (G_\RR, \AA^r_X (\CC), \RR (n)) [-1] \end{array} \ar{d}{\text{split}}[swap]{\cong} \ar{r}{\cong} & \begin{array}{c} \RHom (R\Gamma (X_\et, \ZZ^c (n-r)) [2r], \RR) [-1] \\ \oplus \\ R\Gamma_c (G_\RR, X (\CC), \RR (n-r)) [-1] [-2r] \end{array} \ar{d}{\text{split}}[swap]{\cong} \\
        R\Gamma_\Wc (\AA^r_X, \RR (n)) \ar{r}{\cong} & R\Gamma_\Wc (X, \RR (n-r)) [-2r]
      \end{tikzcd} \]
    Here the first square is commutative due to the compatibility of the
    regulator with affine bundles (Lemma~\ref{lemma:compatibility-of-B(X,n)}),
    and the second square commutes because the quasi-isomorphism
    \eqref{eqn:RGamma-Wc-and-affine-bundles} gives compatible splittings
    (see again Figure~\ref{fig:RGamma-Wc-and-affine-bundles} on
    p.~\pageref{fig:RGamma-Wc-and-affine-bundles}). Taking the determinants, we
    obtain the desired commutative diagram
    \eqref{eqn:lambda-and-affine-bundles}.
  \end{proof}
\end{lemma}

\begin{landscape}
  \begin{figure}
    \[ \begin{tikzcd}[font=\small,column sep=1em]
        &[0.5em] &[-2.75em] &[-2.75em] &[-2.75em] R\Gamma_\Wc (\AA^r_X, \ZZ (n))\ar{dl}\ar[dashed,near start]{dd}{\cong} &[-2.75em] \\
        \RHom (R\Gamma (\AA^r_{X,\et}, \ZZ^c (n)), \QQ[-2]) \ar{r}{\alpha_{\AA^r_X,n}}\ar{dd}{(p^*)^\vee}[swap]{\cong} & R\Gamma_c (\AA^r_{X,\et}, \ZZ (n)) \ar{dr}{u_\infty} \ar{dd}{p_*}[swap]{\cong} \ar{rr} & & R\Gamma_\fg (\AA^r_X, \ZZ (n))\ar{dl}[swap]{i_\infty} \ar[dashed]{dd}{\cong} \ar[crossing over]{rr} & & \cdots [-1] \ar{dd}{\cong} \\
        & & R\Gamma_c (G_\RR, \AA^r_X (\CC), \ZZ (n)) & & R\Gamma_\Wc (X, \ZZ (n-r)) [-2r]\ar{dl} \\
        \RHom (R\Gamma (X_\et, \ZZ^c (n-r)) [2r], \QQ[-2]) \ar[outer sep=0.75em]{r}{\alpha_{X,n-r} [-2r]} & R\Gamma_c (X_\et, \ZZ (n-r)) [-2r] \ar{dr}{u_\infty} \ar{rr} & & R\Gamma_\fg (X, \ZZ (n-r)) [-2r] \ar{dl}[swap]{i_\infty} \ar{rr} & & \cdots [-1] \\
        & & R\Gamma_c (G_\RR, X (\CC), \ZZ (n-r)) [-2r] \ar[<-,near end,crossing over]{uu}[swap]{p_*}{\cong}\\
        & & & & R\Gamma_\Wc (\AA^r_X, \RR (n))\ar{dl}\ar[dashed,near start]{dd}{\cong} & \\
        \RHom (R\Gamma (\AA^r_{X,\et}, \ZZ^c (n)), \RR[-2]) \ar{r}\ar{dd}{(p^*)^\vee}[swap]{\cong} & 0 \ar{dr} \ar{dd} \ar{rr} & & R\Gamma_\fg (\AA^r_X, \RR (n))\ar{dl}[swap]{0} \ar[dashed]{dd}{\cong} \ar[crossing over]{rr} & & \cdots [-1] \ar{dd}{\cong} \\
        & & R\Gamma_c (G_\RR, \AA^r_X (\CC), \RR (n)) & & R\Gamma_\Wc (X, \RR (n-r)) [-2r]\ar{dl} \\
        \RHom (R\Gamma (X_\et, \ZZ^c (n-r)) [2r], \RR[-2]) \ar{r} & 0 \ar{dr} \ar{rr} & & R\Gamma_\fg (X, \RR (n-r)) [-2r] \ar{dl}[swap]{0} \ar{rr} & & \cdots [-1] \\
        & & R\Gamma_c (G_\RR, X (\CC), \RR (n-r)) [-2r] \ar[<-,near end,crossing over]{uu}[swap]{p_*}{\cong}
      \end{tikzcd} \]

    \caption{Isomorphism
      $R\Gamma_\Wc (\AA^r_X, \ZZ (n)) \cong R\Gamma_\Wc (X, \ZZ (n-r)) [-2r]$
      and its splitting after tensoring with $\RR$.}
    \label{fig:RGamma-Wc-and-affine-bundles}
  \end{figure}
\end{landscape}

\begin{theorem}
  \label{thm:compatibility-of-C(X,n)}
  For an arithmetic scheme $X$ and $n < 0$, assume $\mathbf{L}^c (X_\et,n)$,
  $\mathbf{B} (X,n)$, and the meromorphic continuation of $\zeta (X,s)$ around
  $s = n$.

  \begin{enumerate}
  \item[1)] If $X = \coprod_{1 \le i \le r} X_i$ is a finite disjoint union of
    arithmetic schemes, then
    $$\mathbf{C} (X,n) \iff \mathbf{C} (X_i,n)\text{ for all }i.$$

  \item[2)] For a closed-open decomposition
    $Z \not\hookrightarrow X \hookleftarrow U$, if two of three conjectures
    \[ \mathbf{C} (X,n), \quad
      \mathbf{C} (Z,n), \quad
      \mathbf{C} (U,n) \]
    are true, then the third is also true.

  \item[3)] For any $r \ge 0$, one has
    $$\mathbf{C} (\AA^r_X, n) \iff \mathbf{C} (X, n-r).$$
  \end{enumerate}

  \begin{proof}
    Follows from Lemmas
    \ref{lemma:lambda-and-disjoint-unions},
    \ref{lemma:lambda-and-closed-open-decompositions},
    \ref{lemma:lambda-and-affine-bundles},
    together with the corresponding identities for the zeta functions
    \eqref{eqn:zeta-function-for-disjoint-unions},
    \eqref{eqn:zeta-function-for-closed-open-decompositions},
    \eqref{eqn:zeta-function-for-affine-space}.
  \end{proof}
\end{theorem}

The following is a special case of compatibility with closed-open
decompositions.

\begin{lemma}
  \label{lemma:compatibility-for-Xred}
  For an arithmetic scheme $X$ and $n < 0$, the conjectures
  \[
    \mathbf{L}^c (X_\et,n), ~
    \mathbf{B} (X,n), ~
    \mathbf{VO} (X,n), ~
    \mathbf{C} (X,n)
  \]
  are equivalent to
  \[
    \mathbf{L}^c (X_\text{\it red,\'{e}t},n), ~
    \mathbf{B} (X_\red,n), ~
    \mathbf{VO} (X_\red,n), ~
    \mathbf{C} (X_\red,n)
  \]
  respectively.

  \begin{proof}
    Apply Lemma~\ref{lemma:compatibility-of-Lc(X,n)},
    Lemma~\ref{lemma:compatibility-of-B(X,n)},
    Proposition~\ref{prop:compatibility-of-VO(X,n)}, and
    Theorem~\ref{thm:compatibility-of-C(X,n)} to the canonical closed embedding
    $X_\red \hookrightarrow X$.
  \end{proof}
\end{lemma}

The above lemma can be proved directly, by going through the construction of
Weil-\'{e}tale cohomology in \cite{Beshenov-Weil-etale-1} and the statements of
the conjectures. In particular,
$$R\Gamma_\Wc (X, \ZZ(n)) \cong R\Gamma_\Wc (X_\red, \ZZ(n)).$$
It is important to note that the cycle complexes do not distinguish $X$ from
$X_\red$, and neither does the zeta function: $\zeta (X,s) = \zeta (X_\red,s)$.

\begin{remark}
  If $X/\FF_q$ is a variety over a finite field, then the proof of
  Theorem~\ref{thm:compatibility-of-C(X,n)} simplifies drastically: we can work
  with the formula \eqref{eqn:special-value-for-X/Fq} and the following
  properties of motivic cohomology:
  \begin{enumerate}
  \item[1)] $R\Gamma (\coprod_i X_{i,\et}, \ZZ^c (n)) \cong
    \bigoplus_i R\Gamma (X_{i,\et}, \ZZ^c (n))$;

  \item[2)] triangles associated to closed-open decompositions
    \[ R\Gamma (Z_\et, \ZZ^c (n)) \to
      R\Gamma (X_\et, \ZZ^c (n)) \to
      R\Gamma (U_\et, \ZZ^c (n)) \to
      R\Gamma (Z_\et, \ZZ^c (n))[1] \]

  \item[3)] homotopy invariance
    $R\Gamma (X_\et, \ZZ^c (n-r)) [2r] \cong
    R\Gamma (\AA^r_{X,\et}, \ZZ^c (n))$.
  \end{enumerate}
  In this case, no regulators are involved, so we do not need the technical
  lemmas \ref{lemma:lambda-and-disjoint-unions},
  \ref{lemma:lambda-and-closed-open-decompositions},
  \ref{lemma:lambda-and-affine-bundles}.
\end{remark}

If we consider the projective space $\PP_X^r = \PP_\ZZ^r \times X$, we have a
formula for the zeta function
\begin{equation}
  \label{eqn:zeta-function-for-projective-space}
  \zeta (\PP_X^r, s) = \prod_{0 \le i \le r} \zeta (X, s-i).
\end{equation}
Our special value conjecture satisfies the corresponding compatibility.

\begin{corollary}[Projective bundles]
  Let $X$ be an arithmetic scheme, $n < 0$, and $r \ge 0$.
  For all $0 \le i \le r$ assume Conjectures $\mathbf{L}^c (X_\et,n-i)$,
  $\mathbf{B} (X,n-i)$, and the meromorphic continuation of $\zeta (X,s)$ around
  $s = n-i$. Then
  \[ \mathbf{C} (X,n-i)\text{ for }0 \le i \le r \Longrightarrow
    \mathbf{C} (\PP_X^r, n). \]

  \begin{proof}
    Applied to the closed-open decomposition
    $\PP_X^{r-1} \not\hookrightarrow \PP_X^r \hookleftarrow \AA_X^r$,
    Theorem~\ref{thm:compatibility-of-C(X,n)} gives
    \[ \mathbf{C} (X, n-r) \text{ and } \mathbf{C} (\PP_X^{r-1}, n)
      \Longrightarrow
      \mathbf{C} (\AA_X^r, n) \text{ and } \mathbf{C} (\PP_X^{r-1}, n)
      \Longrightarrow
      \mathbf{C} (\PP_X^r,n). \]
    The assertion follows by induction on $r$.
    (The same inductive argument proves the identity
    \eqref{eqn:zeta-function-for-projective-space} from
    \eqref{eqn:zeta-function-for-affine-space}.)
    \end{proof}
\end{corollary}

%%%%%%%%%%%%%%%%%%%%%%%%%%%%%%%%%%%%%%%%%%%%%%%%%%%%%%%%%%%%%%%%%%%%%%%%%%%%%%%%

\section{Unconditional results}
\label{sec:unconditional-results}

Now we apply Theorem~\ref{thm:compatibility-of-C(X,n)} to prove the main theorem
stated in the introduction: the validity of $\mathbf{VO} (X,n)$ and
$\mathbf{C} (X,n)$ for all $n < 0$ for cellular schemes over certain
one-dimensional bases. In fact, we will construct an even larger class of
schemes $\mathcal{C} (\ZZ)$ whose elements satisfy the conjectures. This
approach is motivated by \cite[\S 5]{Morin-2014}.

\begin{definition}
  Let $\mathcal{C} (\ZZ)$ be the full subcategory of the category of arithmetic
  schemes generated by the following objects:
  \begin{itemize}
  \item the empty scheme $\emptyset$,
  \item $\Spec \FF_q$ for each finite field,
  \item $\Spec \mathcal{O}_F$ for an abelian number field $F/\QQ$,
  \item curves over finite fields $C/\FF_q$,
  \end{itemize}
  and the following operations.
  \begin{enumerate}
  \item[$\mathcal{C}0)$] $X$ is in $\mathcal{C} (\ZZ)$ if and only if $X_\red$
    is in $\mathcal{C} (\ZZ)$.

  \item[$\mathcal{C}1)$] A finite disjoint union $\coprod_{1 \le i \le r} X_i$
    is in $\mathcal{C} (\ZZ)$ if and only if each $X_i$ is in
    $\mathcal{C} (\ZZ)$.

  \item[$\mathcal{C}2)$] Let $Z \not\hookrightarrow X \hookleftarrow U$ be a
    closed-open decomposition such that $Z_{\red,\CC}$, $X_{\red,\CC}$,
    $U_{\red,\CC}$ are smooth and quasi-projective. If two of the three schemes
    $Z,X,U$ lie in $\mathcal{C} (\ZZ)$, then the third also lies in
    $\mathcal{C} (\ZZ)$.

  \item[$\mathcal{C}3)$] If $X$ lies in $\mathcal{C} (\ZZ)$, then the affine
    space $\AA^r_X$ for each $r \ge 0$ also lies in $\mathcal{C} (\ZZ)$.
  \end{enumerate}
\end{definition}

Recall that the condition that $X_{\red,\CC}$ is smooth and quasi-projective is
necessary to ensure that the regulator morphism exists (see
Remark~\ref{rmk:regulator-is-defined-for-XC-smooth-quasi-proj}).

\begin{proposition}
  \label{prop:C(X,n)-holds-for-C(Z)}
  Conjectures $\mathbf{VO} (X,n)$ and $\mathbf{C} (X,n)$ hold for any
  $X \in \mathcal{C} (\ZZ)$ and $n < 0$.

  \begin{proof}
    Finite fields satisfy $\mathbf{C} (X,n)$ by
    Example~\ref{example:C(X,n)-for-Spec-Fq}.

    If $X = \Spec \mathcal{O}_F$ for a number field $F/\QQ$, then
    Conjecture~$\mathbf{C} (X,n)$ is equivalent to the conjecture of Flach and
    Morin \cite[Conjecture~5.12]{Flach-Morin-2018}, which holds unconditionally
    for abelian $F/\QQ$, via reduction to the Tamagawa number conjecture; see
    \cite[\S 5.8.3]{Flach-Morin-2018}, in particular
    [ibid., Proposition~5.35]. The condition $\mathbf{VO} (X,n)$ is also true in
    this case (see Example~\ref{example:VO(X,n)-for-number-rings}).

    If $X = C/\FF_q$ is a curve over a finite field, then $\mathbf{C} (X,n)$
    holds thanks to
    Theorem~\ref{thm:C(X,n)-over-finite-fields}.
    Conjecture~$\mathbf{L}^c (X_\et,n)$ is known for curves and essentially goes
    back to Soul\'{e}; see, for example, \cite[Proposition~4.3]{Geisser-2017}.

    Finally, the fact that the Conjectures $\mathbf{L}^c (X_\et,n)$,
    $\mathbf{B} (X,n)$, $\mathbf{VO} (X,n)$, $\mathbf{C} (X,n)$ are closed under
    the operations $\mathcal{C}0)$--$\mathcal{C}3)$ is
    Lemma~\ref{lemma:compatibility-for-Xred},
    Lemma~\ref{lemma:compatibility-of-Lc(X,n)},
    Lemma~\ref{lemma:compatibility-of-B(X,n)},
    Proposition~\ref{prop:compatibility-of-VO(X,n)}, and
    Theorem~\ref{thm:compatibility-of-C(X,n)} respectively.
  \end{proof}
\end{proposition}

\begin{lemma}
  Any zero-dimensional arithmetic scheme $X$ is in $\mathcal{C} (\ZZ)$.

  \begin{proof}
    Since $X$ is a Noetherian scheme of dimension $0$, it is a finite disjoint
    union of $\Spec A_i$ for some Artinian local rings $A_i$. Thanks to
    $\mathcal{C}1)$, we can assume that $X = \Spec A$, and thanks to
    $\mathcal{C}0)$, we can assume that $X$ is reduced. But then $A = k$ is a
    field. Since $X$ is a scheme of finite type over $\Spec \ZZ$, we conclude
    that $X = \Spec \FF_q \in \mathcal{C} (\ZZ)$.
  \end{proof}
\end{lemma}

\begin{proposition}
  \label{prop:particular-cases-1-dim-base}
  Let $B$ be a one-dimensional arithmetic scheme. Suppose that each of the
  generic points $\eta \in B$ satisfies one of the following properties:
  \begin{enumerate}
  \item[a)] $\fchar \kappa (\eta) = p > 0$;

  \item[b)] $\fchar \kappa (\eta) = 0$, and $\kappa (\eta)/\QQ$ is an abelian
    number field.
  \end{enumerate}
  Then $B \in \mathcal{C} (\ZZ)$.

  \begin{proof}
    We verify that such a scheme can be obtained from $\Spec \mathcal{O}_F$ for
    an abelian number field $F/\QQ$, or a curve over a finite field $C/\FF_q$,
    using the operations $\mathcal{C}0)$, $\mathcal{C}1)$, $\mathcal{C}2)$ which
    appear in the definition of $\mathcal{C} (\ZZ)$.

    Thanks to $\mathcal{C}0)$, we can assume that $B$ is reduced. Consider the
    normalization $\nu\colon B' \to B$. This is a birational morphism, so there
    exist open dense subsets $U' \subseteq B'$ and $U \subseteq B$ such that
    $\left.\nu\right|_{U'}\colon U' \xrightarrow{\cong} U$. Now $B\setminus U$
    is zero-dimensional, and therefore $B\setminus U \in \mathcal{C} (\ZZ)$ by
    the previous lemma. Thanks to $\mathcal{C}2)$, it suffices to check that
    $U' \in \mathcal{C} (\ZZ)$, and this would imply $B \in \mathcal{C} (\ZZ)$.

    Now $U'$ is a finite disjoint union of normal integral schemes, so according
    to $\mathcal{C}1)$ we can assume that $U'$ is integral. Consider the generic
    point $\eta \in U'$ and the residue field $F = \kappa (\eta)$. There are two
    cases to consider.

    \begin{enumerate}
    \item[a)] If $\fchar F = p > 0$, then $U'$ is a curve over a finite field,
      so it lies in $\mathcal{C} (\ZZ)$.

    \item[b)] If $\fchar F = 0$, then by our assumptions, $F/\QQ$ is an abelian
      number field.

      We note that if $V' \subseteq U'$ is an affine open neighborhood of
      $\eta$, then $U'\setminus V' \in \mathcal{C} (\ZZ)$ by the previous
      lemma. Therefore, we can assume without loss of generality that $U'$ is
      affine.

      We have $U' = \Spec \mathcal{O}$, where $\mathcal{O}$ is a finitely
      generated integrally closed domain. This means that
      $\mathcal{O}_F \subseteq \mathcal{O} = \mathcal{O}_{F,S}$ for a finite set
      of places $S$. Now $U' = \Spec \mathcal{O}_F \setminus S$, and
      $S \in \mathcal{C} (\ZZ)$, so everything reduces to the case of
      $U' = \Spec \mathcal{O}_F$, which is in $\mathcal{C} (\ZZ)$. \qedhere
    \end{enumerate}
  \end{proof}
\end{proposition}

\begin{remark}
  Schemes like the above were considered by Jordan and Poonen in
  \cite{Jordan-Poonen-2020}, where the authors write down a special value
  formula for $s = 1$ that generalizes the classical class number
  formula. Namely, they consider the case where $B$ is reduced and affine,
  but without requiring $\kappa (\eta)/\QQ$ to be abelian.
\end{remark}

\begin{example}
  If $B = \Spec \mathcal{O}$ for a nonmaximal order
  $\mathcal{O} \subset \mathcal{O}_F$, where $F/\QQ$ is an abelian number field,
  then our formalism gives a cohomological interpretation of the special values
  of $\zeta_\mathcal{O} (s)$ at $s = n < 0$. This already seems to be a new
  result.
\end{example}

\begin{definition}
  \label{dfn:B-cellular-scheme}
  Let $X \to B$ be a $B$-scheme. We say that $X$ is \textbf{$B$-cellular} if it
  admits a filtration by closed subschemes
  \begin{equation}
    \label{eqn:cellular-decomposition}
    X = Z_N \supseteq Z_{N-1} \supseteq \cdots \supseteq Z_0 \supseteq Z_{-1} = \emptyset
  \end{equation}
  such that $Z_i\setminus Z_{i-1} \cong \coprod_j \AA^{r_{i_j}}_B$ is a finite
  union of affine $B$-spaces.
\end{definition}

For example, projective spaces $\PP^r_B$ and, in general, Grassmannians
$\Gr (k,\ell)_B$ are cellular. Many interesting examples of cellular schemes
arise from actions of algebraic groups on varieties and the Bia\l{}ynicki-Birula
theorem; see \cite{Wendt-2010} and \cite{Brosnan-2005}.

\begin{proposition}
  \label{prop:cellular-schemes-in-C(Z)}
  Let $X$ be a $B$-cellular arithmetic scheme, where $B \in \mathcal{C} (\ZZ)$,
  and $X_{\red,\CC}$ is smooth and quasi-projective. Then
  $X \in \mathcal{C} (\ZZ)$.

  \begin{proof}
    Considering the corresponding cellular decomposition
    \eqref{eqn:cellular-decomposition}, we pass to open complements
    $U_i = X\setminus Z_i$ to obtain a filtration
    $$X = U_{-1} \supseteq U_1 \supseteq \cdots \supseteq U_{N-1} \supseteq U_N = \emptyset,$$
    where $U_{i,\CC}$ are smooth and quasi-projective, being \emph{open}
    subvarieties in $X_\CC$. Now we have closed-open decompositions
    $\coprod_j \AA^{r_{i,j}}_B \not\hookrightarrow U_i \hookleftarrow U_{i+1}$,
    and the claim follows by induction on the length of the cellular
    decomposition, using operations $\mathcal{C}1)$--$\mathcal{C}3)$.
  \end{proof}
\end{proposition}

As a corollary of the above, we obtain the following result, stated in the
introduction.

\begin{theorem}
  Let $B$ be a one-dimensional arithmetic scheme satisfying the assumptions of
  Proposition~\ref{prop:particular-cases-1-dim-base}. If $X$ is a $B$-cellular
  arithmetic scheme with smooth and quasi-projective fiber $X_{\red,\CC}$, then
  Conjectures $\mathbf{VO} (X,n)$ and $\mathbf{C} (X,n)$ hold unconditionally
  for any $n < 0$.

  \begin{proof}
    Follows from propositions
    \ref{prop:C(X,n)-holds-for-C(Z)},
    \ref{prop:particular-cases-1-dim-base},
    \ref{prop:cellular-schemes-in-C(Z)}.
  \end{proof}
\end{theorem}

%%%%%%%%%%%%%%%%%%%%%%%%%%%%%%%%%%%%%%%%%%%%%%%%%%%%%%%%%%%%%%%%%%%%%%%%%%%%%%%%

\begin{appendices}
\section{Determinants of complexes}
\label{app:determinants}

Here we give a brief overview of the determinants of complexes. The original
construction goes back to Knudsen and Mumford \cite{Knudsen-Mumford-1976}, and
useful expositions can be found in
\cite[Appendix~A]{Gelfand-Kapranov-Zelevinsky-1994} and \cite[\S
2.1]{Kato-1993}.

For our purposes, let $R$ be an integral domain.

\begin{definition}
  Denote by $\mathcal{P}_\is (R)$ the category of
  \textbf{graded invertible $R$-modules}. It has as objects $(L,r)$, where $L$
  is an invertible $R$-module (i.e. projective of rank $1$) and $r \in \ZZ$. The
  morphisms in this category are given by
  \[ \Hom_{\mathcal{P}_\is (R)} ((L,r), (M,s)) = \begin{cases}
      \Isom_R (L,M), & r = s, \\
      \emptyset, & r \ne s.
    \end{cases} \]
  This category is equipped with tensor products
  $$(L,r) \otimes_R (M,s) = (L\otimes_R M, r + s)$$
  with (graded) commutativity isomorphisms
  \[
    (L,r) \otimes_R (M,s)
    \xrightarrow{\cong}
    (M,s) \otimes_R (L,r),
    \quad
    \ell \otimes m \mapsto (-1)^{rs}\,m\otimes \ell.
  \]
\end{definition}

The unit object with respect to this product is $(R,0)$, and for each
$(L,r) \in \mathcal{P}_\is (R)$ the inverse is given by $(L^{-1}, -r)$ where
$L^{-1} = \iHom_R (L,R)$. The canonical evaluation morphism
$L \otimes_R \iHom_R (L,R) \to R$ induces an isomorphism
$$(L,r) \otimes_R (L^{-1}, -r) \cong (R,0).$$

\begin{definition}
  Denote by $\mathcal{C}_\is (R)$ the category whose objects are finitely
  generated projective $R$-modules and whose morphisms are isomorphisms.
  For $A \in \mathcal{C}_\is (R)$ we define the corresponding determinant by
  \begin{equation}
    \label{eqn:det-proj-module}
    \det_R (A) = \Bigl(\bigwedge^{\rk_R A}_R A, \rk_R A\Bigr) \in \mathcal{P}_\is (R).
  \end{equation}
  Here $\rk_R A$ is the rank of $A$, so that the top exterior power
  $\bigwedge^{\rk_R A}_R A$ is an invertible $R$-module.
\end{definition}

This yields a functor
$\det_R\colon \mathcal{C}_\is (R) \to \mathcal{P}_\is (R)$.
For $(L,r) \in \mathcal{P}_\is (R)$ we usually forget about $r$ and treat the
determinant as an invertible $R$-module.

The main result of \cite[Chapter~I]{Knudsen-Mumford-1976} is that this
construction can be generalized to complexes and morphisms in the derived
category.

\begin{definition}
  Let $\mathbf{D} (R)$ be the derived category of the category of
  $R$-modules. Recall that a complex $A^\bullet$ is \textbf{perfect} if it is
  quasi-isomorphic to a bounded complex of finitely generated projective
  $R$-modules. We denote by $\Parf_\is (R)$ the subcategory of $\mathbf{D} (R)$
  whose objects consist of perfect complexes, and whose morphisms are
  quasi-isomorphisms of complexes.
\end{definition}

\begin{theorem}[Knudsen--Mumford]
  The determinant \eqref{eqn:det-proj-module} extends to perfect complexes of
  $R$-modules as a functor
  $$\det_R\colon \Parf_\is (R) \to \mathcal{P}_\is (R),$$
  satisfying the following properties.

  \begin{itemize}
  \item $\det_R (0) = (R,0)$.

  \item For a distinguished triangle of complexes in $\Parf_\is (R)$
    \[ A^\bullet \xrightarrow{u}
      B^\bullet \xrightarrow{v}
      C^\bullet \xrightarrow{w} A^\bullet [1] \]
    there is a canonical isomorphism
    \[ i_R (u,v,w)\colon \det_R A^\bullet \otimes_R \det_R C^\bullet
      \xrightarrow{\cong} \det_R B^\bullet. \]

  \item In particular, there exist canonical isomorphisms
    \[ \det_R (A^\bullet \oplus B^\bullet) \cong
      \det_R (A^\bullet) \otimes_R \det_R (B^\bullet). \]

  \item For the triangles
    \[ \begin{tikzcd}[column sep=1em]
        A^\bullet \ar{r}{id} & A^\bullet \ar{r} & 0^\bullet \ar{r} & A^\bullet [1] \\[-2em]
        0^\bullet \ar{r} & A^\bullet \ar{r}{id} & A^\bullet \ar{r} & 0^\bullet [1]
      \end{tikzcd} \]
    the isomorphism $i_R$ comes from the canonical isomorphism
    $\det_R A^\bullet \otimes_R (R,0) \cong \det_R A^\bullet$.

  \item For an isomorphism of distinguished triangles
    \[ \begin{tikzcd}
        A^\bullet \ar{r}{u}\ar{d}{f}[swap]{\cong} & B^\bullet \ar{r}{v}\ar{d}{g}[swap]{\cong} & C^\bullet \ar{r}{w}\ar{d}{h}[swap]{\cong} & A^\bullet [1]\ar{d}{f[1]}[swap]{\cong} \\
        A'^\bullet \ar{r}{u'} & B'^\bullet \ar{r}{v'} & C'^\bullet \ar{r}{w'} & A'^\bullet [1]
      \end{tikzcd} \]
    the diagram
    \[ \begin{tikzcd}[column sep=5em]
        \det_R A^\bullet \otimes_R \det_R C^\bullet \ar{r}{i_R (u,v,w)}[swap]{\cong}\ar{d}{\det_R (f) \otimes \det_R (h)}[swap]{\cong} & \det_R B^\bullet\ar{d}{\det_R (g)}[swap]{\cong} \\
        \det_R A'^\bullet \otimes_R \det_R C'^\bullet \ar{r}{i_R (u',v',w')}[swap]{\cong} & \det_R B'^\bullet
      \end{tikzcd} \]
    is commutative.

  \item The determinant is compatible with base change: given a ring
    homomorphism $R\to S$, there is a natural isomorphism
    \[ \det_S (A^\bullet \otimes_R^\mathbf{L} S)
      \xrightarrow{\cong}
      (\det_R A^\bullet) \otimes_R S. \]
    Moreover, this isomorphism is compatible with $i_R$ and $i_S$.

  \item If $A^\bullet$ is a bounded complex where each object $A^i$ is perfect
    (i.e. admits a finite length resolution by finitely generated projective
    $R$-modules), then
    $$\det_R A^\bullet \cong \bigotimes_{i\in \ZZ} (\det_R A^i)^{(-1)^i}.$$
    If each $A^i$ is already a finitely generated projective $R$-module, then
    $\det_R A^i$ in the above formula is given by \eqref{eqn:det-proj-module}.

  \item If the cohomology modules $H^i (A^\bullet)$ are perfect, then
    \begin{equation}
      \label{eqn:det-in-terms-of-cohomology}
      \det_R A^\bullet \cong
      \bigotimes_{i\in \ZZ} (\det_R H^i (A^\bullet))^{(-1)^i}.
    \end{equation}
  \end{itemize}
\end{theorem}

We refer the reader to \cite{Knudsen-Mumford-1976} for the actual construction
and proofs.

A particularly simple case of interest is when $R = \ZZ$ and all cohomology
groups $H^i (A^\bullet)$ are finite.

\begin{lemma}
  \label{lemma:determinant-for-torsion-cohomology}
  ~

  \begin{enumerate}
  \item[1)] Let $A$ be a finite abelian group. Then
    \[ (\det_\ZZ A) \subset (\det_\ZZ A) \otimes \QQ
      \cong \det_\QQ (A \otimes \QQ) = \det_\QQ (0) \cong \QQ \]
    corresponds to the fractional ideal $\frac{1}{\# A} \ZZ \subset \QQ$.

  \item[2)] In general, let $A^\bullet$ be a perfect complex of abelian groups
    such that the cohomology groups $H^i (A^\bullet)$ are all finite. Then
    $\det_\ZZ A^\bullet$ corresponds to the fractional ideal
    $\frac{1}{m}\,\ZZ \subset \QQ$, where
    $$m = \prod_{i\in \ZZ} |H^i (A^\bullet)|^{(-1)^i}.$$
  \end{enumerate}

  \begin{proof}
    Since $\det_\ZZ (A\oplus B) \cong \det_\ZZ A \otimes \det_\ZZ B$, in
    part 1) it suffices to consider the case of a cyclic group
    $A = \ZZ/m\ZZ$. Using the resolution
    \[ \ZZ/m\ZZ [0] \cong \Bigl[
      \mathop{m\ZZ}_{\text{deg.\,}-1} \hookrightarrow
      \mathop{\ZZ}_{\text{deg.\,}0}
      \Bigr], \]
    we calculate
    $$\det_\ZZ (\ZZ/m\ZZ) \cong \ZZ \otimes (m\ZZ)^{-1} \cong (m\ZZ)^{-1},$$
    which corresponds to $\frac{1}{m}\,\ZZ$ in $\QQ$.
    Part 2) follows directly from 1) and \eqref{eqn:det-in-terms-of-cohomology}.
  \end{proof}
\end{lemma}

\end{appendices}

%%%%%%%%%%%%%%%%%%%%%%%%%%%%%%%%%%%%%%%%%%%%%%%%%%%%%%%%%%%%%%%%%%%%%%%%%%%%%%%%

\bibliographystyle{abbrv}
\bibliography{../weil-etale}

\begin{thebibliography}{10}

\bibitem{Bayer-Neukirch-1978}
P.~B\'{a}yer and J.~Neukirch.
\newblock On values of zeta functions and {$l$}-adic {E}uler characteristics.
\newblock {\em Invent. Math.}, 50(1):35--64, 1978.

\bibitem{Beshenov-Weil-etale-1}
A.~Beshenov.
\newblock Weil-\'{e}tale cohomology and duality for arithmetic schemes in
  negative weights, 2021.
\newblock preprint.

\bibitem{Bhatt-Scholze-2015}
B.~Bhatt and P.~Scholze.
\newblock The pro-\'{e}tale topology for schemes.
\newblock {\em Ast\'{e}risque}, (369):99--201, 2015.

\bibitem{Bloch-1986}
S.~Bloch.
\newblock Algebraic cycles and higher {$K$}-theory.
\newblock {\em Adv. in Math.}, 61(3):267--304, 1986.

\bibitem{Bloch-1986-Lefschetz}
S.~Bloch.
\newblock Algebraic cycles and the {B}eilinson conjectures.
\newblock In {\em The {L}efschetz centennial conference, {P}art {I} ({M}exico
  {C}ity, 1984)}, volume~58 of {\em Contemp. Math.}, pages 65--79. Amer. Math.
  Soc., Providence, RI, 1986.

\bibitem{Brosnan-2005}
P.~Brosnan.
\newblock On motivic decompositions arising from the method of
  {B}ia\l{}ynicki-{B}irula.
\newblock {\em Invent. Math.}, 161:91--111, 2005.

\bibitem{SGA4-1-2}
P.~Deligne.
\newblock {\em Cohomologie \'etale}.
\newblock Lecture Notes in Mathematics, Vol. 569. Springer-Verlag, Berlin-New
  York, 1977.
\newblock S{\'e}minaire de G{\'e}om{\'e}trie Alg{\'e}brique du {B}ois-{M}arie
  {SGA}~4$\frac{1}{2}$, Avec la collaboration de {J}.~{F}.~{B}outot,
  {A}.~{G}rothendieck, {L}.~Illusie et {J}.~{L}.~{V}erdier.

\bibitem{Deligne-Weil-II}
P.~Deligne.
\newblock La conjecture de {W}eil. {II}.
\newblock {\em Inst. Hautes \'{E}tudes Sci. Publ. Math.}, (52):137--252, 1980.

\bibitem{Esnault-Viehweg-1988}
H.~Esnault and E.~Viehweg.
\newblock Deligne-{B}eilinson cohomology.
\newblock In {\em Beilinson's conjectures on special values of
  {$L$}-functions}, volume~4 of {\em Perspect. Math.}, pages 43--91. Academic
  Press, Boston, MA, 1988.

\bibitem{Flach-Morin-2012}
M.~Flach and B.~Morin.
\newblock On the {W}eil-\'{e}tale topos of regular arithmetic schemes.
\newblock {\em Doc. Math.}, 17:313--399, 2012.

\bibitem{Flach-Morin-2018}
M.~Flach and B.~Morin.
\newblock Weil-\'{e}tale cohomology and zeta-values of proper regular
  arithmetic schemes.
\newblock {\em Doc. Math.}, 23:1425--1560, 2018.

\bibitem{Flach-Morin-2020}
M.~Flach and B.~Morin.
\newblock Compatibility of special value conjectures with the functional
  equation of zeta functions, 2020.
\newblock preprint.

\bibitem{Fontaine-Perrin-Riou-1994}
J.-M. Fontaine and B.~Perrin-Riou.
\newblock Autour des conjectures de {B}loch et {K}ato: cohomologie galoisienne
  et valeurs de fonctions {$L$}.
\newblock In {\em Motives ({S}eattle, {WA}, 1991)}, volume~55 of {\em Proc.
  Sympos. Pure Math.}, pages 599--706. Amer. Math. Soc., Providence, RI, 1994.

\bibitem{Geisser-2004-Dedekind}
T.~Geisser.
\newblock Motivic cohomology over {D}edekind rings.
\newblock {\em Math. Z.}, 248(4):773--794, 2004.

\bibitem{Geisser-2004}
T.~Geisser.
\newblock Weil-\'{e}tale cohomology over finite fields.
\newblock {\em Math. Ann.}, 330(4):665--692, 2004.

\bibitem{Geisser-2006}
T.~Geisser.
\newblock Arithmetic cohomology over finite fields and special values of
  {$\zeta$}-functions.
\newblock {\em Duke Math. J.}, 133(1):27--57, 2006.

\bibitem{Geisser-2010-arithmetic-homology}
T.~Geisser.
\newblock Arithmetic homology and an integral version of {K}ato's conjecture.
\newblock {\em J. Reine Angew. Math.}, 644:1--22, 2010.

\bibitem{Geisser-2010}
T.~Geisser.
\newblock Duality via cycle complexes.
\newblock {\em Ann. of Math. (2)}, 172(2):1095--1126, 2010.

\bibitem{Geisser-2017}
T.~Geisser.
\newblock On the structure of \'{e}tale motivic cohomology.
\newblock {\em J. Pure Appl. Algebra}, 221(7):1614--1628, 2017.

\bibitem{Gelfand-Kapranov-Zelevinsky-1994}
I.~Gel${}'$fand, M.~Kapranov, and A.~Zelevinsky.
\newblock {\em Discriminants, resultants, and multidimensional determinants}.
\newblock Mathematics: Theory \& Applications. Birkh\"{a}user Boston, Inc.,
  Boston, MA, 1994.

\bibitem{Grothendieck-FL}
A.~Grothendieck.
\newblock Formule de {L}efschetz et rationalit\'{e} des fonctions {$L$}.
\newblock In {\em Dix expos\'{e}s sur la cohomologie des sch\'{e}mas}, volume~3
  of {\em Adv. Stud. Pure Math.}, pages 31--45. North-Holland, Amsterdam, 1968.

\bibitem{Jannsen-1988}
U.~Jannsen.
\newblock Deligne homology, {H}odge-{$\mathcal{D}$}-conjecture, and motives.
\newblock In {\em Beilinson's conjectures on special values of
  {$L$}-functions}, volume~4 of {\em Perspect. Math.}, pages 305--372. Academic
  Press, Boston, MA, 1988.

\bibitem{Jordan-Poonen-2020}
B.~W. Jordan and B.~Poonen.
\newblock The analytic class number formula for 1-dimensional affine schemes.
\newblock {\em Bull. Lond. Math. Soc.}, 52(5):793--806, 2020.

\bibitem{Kahn-2005}
B.~Kahn.
\newblock Algebraic {$K$}-theory, algebraic cycles and arithmetic geometry.
\newblock In {\em Handbook of {$K$}-theory. {V}ol. 1, 2}, pages 351--428.
  Springer, Berlin, 2005.

\bibitem{Kato-1993}
K.~Kato.
\newblock Lectures on the approach to {I}wasawa theory for {H}asse-{W}eil
  {$L$}-functions via {$B_{\rm dR}$}. {I}.
\newblock In {\em Arithmetic algebraic geometry ({T}rento, 1991)}, volume 1553
  of {\em Lecture Notes in Math.}, pages 50--163. Springer, Berlin, 1993.

\bibitem{Katz-1994}
N.~M. Katz.
\newblock Review of {$\ell$}-adic cohomology.
\newblock In {\em Motives ({S}eattle, {WA}, 1991)}, volume~55 of {\em Proc.
  Sympos. Pure Math.}, pages 21--30. Amer. Math. Soc., Providence, RI, 1994.

\bibitem{Kedlaya-2006}
K.~S. Kedlaya.
\newblock Fourier transforms and {$p$}-adic `{W}eil {II}'.
\newblock {\em Compos. Math.}, 142(6):1426--1450, 2006.

\bibitem{Kerr-Lewis-2007}
M.~Kerr and J.~D. Lewis.
\newblock The {A}bel-{J}acobi map for higher {C}how groups. {II}.
\newblock {\em Invent. Math.}, 170(2):355--420, 2007.

\bibitem{Kerr-Lewis-Lopatto-2018}
M.~Kerr, J.~D. Lewis, and P.~Lopatto.
\newblock Simplicial {A}bel-{J}acobi maps and reciprocity laws.
\newblock {\em J. Algebraic Geom.}, 27(1):121--172, 2018.
\newblock With an appendix by Jos\'{e} Ignacio Burgos-Gil.

\bibitem{Kerr-Lewis-Muller-Stach-2006}
M.~Kerr, J.~D. Lewis, and S.~M\"{u}ller-Stach.
\newblock The {A}bel-{J}acobi map for higher {C}how groups.
\newblock {\em Compos. Math.}, 142(2):374--396, 2006.

\bibitem{Knudsen-Mumford-1976}
F.~F. Knudsen and D.~Mumford.
\newblock The projectivity of the moduli space of stable curves. {I}.
  {P}reliminaries on ``det'' and ``{D}iv''.
\newblock {\em Math. Scand.}, 39(1):19--55, 1976.

\bibitem{Levine-1994}
M.~Levine.
\newblock Bloch's higher {C}how groups revisited.
\newblock In {\em {$K$}-theory ({S}trasbourg, 1992)}, number 226, pages
  235--320. Soci\'et\'e math\'ematique de France, 1994.

\bibitem{Lichtenbaum-1984}
S.~Lichtenbaum.
\newblock Values of zeta-functions at nonnegative integers.
\newblock In {\em Number theory, {N}oordwijkerhout 1983 ({N}oordwijkerhout,
  1983)}, volume 1068 of {\em Lecture Notes in Math.}, pages 127--138.
  Springer, Berlin, 1984.

\bibitem{Lichtenbaum-2005}
S.~Lichtenbaum.
\newblock The {W}eil-\'{e}tale topology on schemes over finite fields.
\newblock {\em Compos. Math.}, 141(3):689--702, 2005.

\bibitem{Lichtenbaum-2009-Euler-char}
S.~Lichtenbaum.
\newblock Euler characteristics and special values of zeta-functions.
\newblock In {\em Motives and algebraic cycles}, volume~56 of {\em Fields Inst.
  Commun.}, pages 249--255. Amer. Math. Soc., Providence, RI, 2009.

\bibitem{Lichtenbaum-2009-number-rings}
S.~Lichtenbaum.
\newblock The {W}eil-\'{e}tale topology for number rings.
\newblock {\em Ann. of Math. (2)}, 170(2):657--683, 2009.

\bibitem{Lichtenbaum-2021}
S.~Lichtenbaum.
\newblock Special values of zeta functions of schemes, 2021.
\newblock preprint.

\bibitem{Milne-1986}
J.~S. Milne.
\newblock Values of zeta functions of varieties over finite fields.
\newblock {\em Amer. J. Math.}, 108(2):297--360, 1986.

\bibitem{Morin-2014}
B.~Morin.
\newblock Zeta functions of regular arithmetic schemes at {$s=0$}.
\newblock {\em Duke Math. J.}, 163(7):1263--1336, 2014.

\bibitem{Neeman-1991}
A.~Neeman.
\newblock Some new axioms for triangulated categories.
\newblock {\em J. Algebra}, 139(1):221--255, 1991.

\bibitem{Ramachandran-2016}
N.~Ramachandran.
\newblock Higher {E}uler characteristics: variations on a theme of {E}uler.
\newblock {\em Homology Homotopy Appl.}, 18(1):231--246, 2016.

\bibitem{Serre-1965}
J.-P. Serre.
\newblock Zeta and {$L$} functions.
\newblock In {\em Arithmetical {A}lgebraic {G}eometry ({P}roc. {C}onf.
  {P}urdue}, pages 82--92. Harper \& Row, New York, 1965.

\bibitem{Serre-1970}
J.-P. Serre.
\newblock Facteurs locaux des fonctions z\^eta des vari\'et\'es alg\'ebriques
  (d\'efinitions et conjectures).
\newblock {\em S\'eminaire Delange-Pisot-Poitou. Th\'eorie des nombres}, 11(2),
  1969--1970.

\bibitem{Siebel-2019}
D.~A. Siebel.
\newblock Special values of zeta-functions for proper regular arithmetic
  surfaces, 2019.
\newblock PhD thesis, California Institute of Technology.

\bibitem{Soule-1984-ICM}
C.~Soul\'{e}.
\newblock {$K$}-th\'{e}orie et z\'{e}ros aux points entiers de fonctions
  z\^{e}ta.
\newblock In {\em Proceedings of the {I}nternational {C}ongress of
  {M}athematicians, {V}ol. 1, 2 ({W}arsaw, 1983)}, pages 437--445. PWN, Warsaw,
  1984.

\bibitem{Stacks-project}
{The Stacks Project Authors}.
\newblock Stacks project, 2021.

\bibitem{Weisschuh-2017}
T.~Wei\ss{}schuh.
\newblock A commutative regulator map into {D}eligne-{B}eilinson cohomology.
\newblock {\em Manuscripta Math.}, 152(3-4):281--315, 2017.

\bibitem{Wendt-2010}
M.~Wendt.
\newblock More examples of motivic cell structures, 2010.
\newblock preprint.

\end{thebibliography}

\end{document}